\documentclass[english,svgnames,table]{article}
\usepackage{preamble}

\declaretheorem[style=definition, name=Theorem, refname={Theorem,Theorems}]{alphThmPrime}

\title{Representability of Galois extensions}

\begin{document}
\begin{titlepage}
    \maketitle
    \thispagestyle{empty}
    \begin{abstract}
        We study faithful Galois extensions in a presentably symmetric monoidal $\infty$-category $\cC$ by a dualizable $\EE_1$-group $G$.
        We prove that these extensions correspond to colimit-preserving symmetric monoidal functors from $\An^{\BG}$ to $\cC$.
        Under suitable assumptions, we use this correspondence to construct a Galois pro-anima that classifies Galois extensions in $\cC$.
        This generalizes the \'etale fundamental groupoid of commutative rings and Mathew's weak Galois groupoid.
        \begin{figure}[h]
            \centering
            \includegraphics[width=0.8\linewidth]{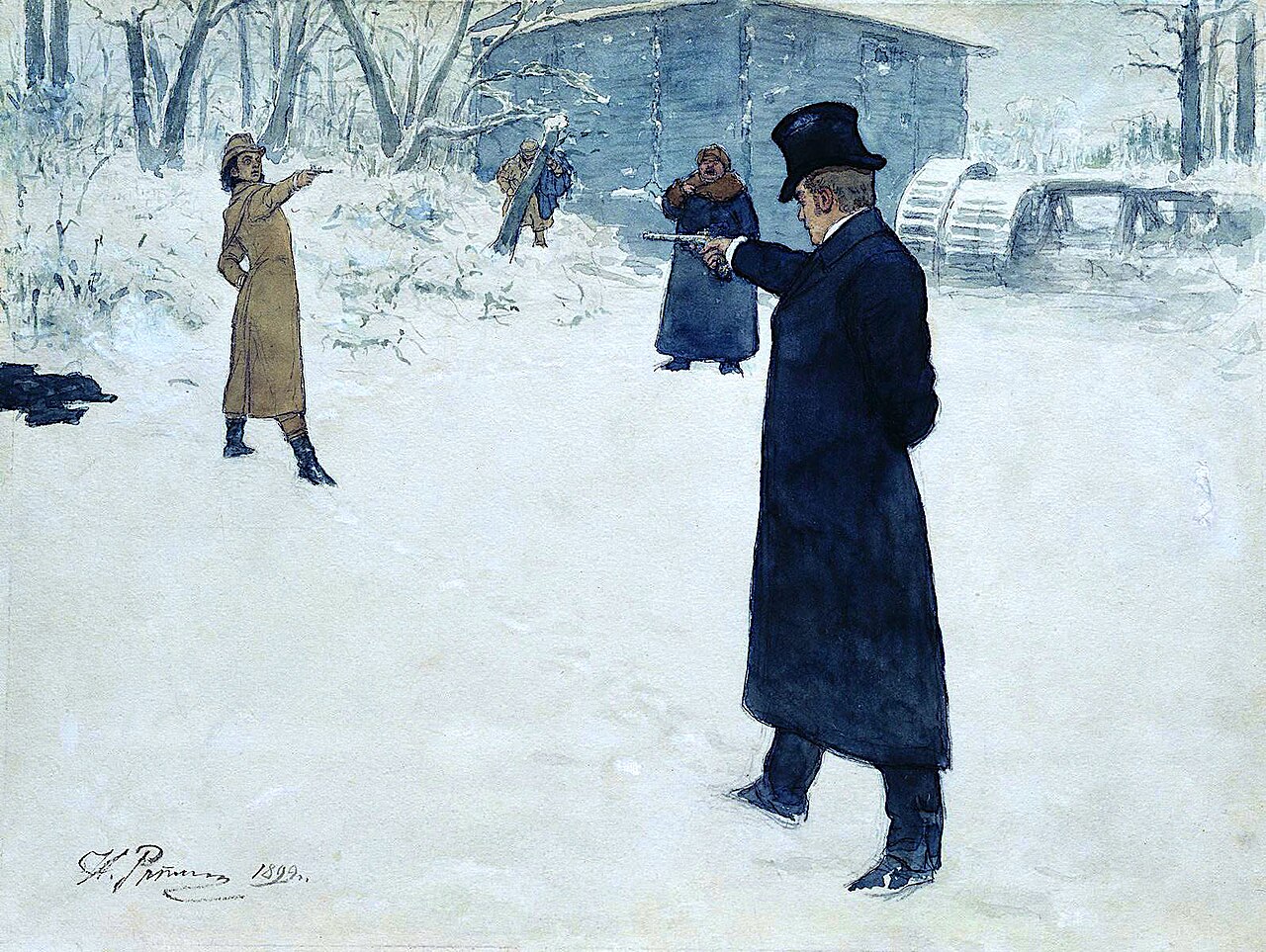}
            \caption*{Ilya Repin, \textit{Duel}, 1899}
        \end{figure}
    \end{abstract}    
\end{titlepage}

\tableofcontents

\section{Introduction}
\label{sec:intro}
Classical Galois theory studies a field $K$ through its finite Galois extensions and their symmetries. 
These extensions are organized by the profinite absolute Galois group $\G_K$. Extensions with a Galois group $G$, for $G$ finite, are classified by continuous surjective homomorphisms $\G_K \onto G$, up to conjugation.

Finite Galois theory was extended to commutative rings in~\cite{Auslander-Goldman-1960-Brauer,Chase-Harrison-Rosenberg-1965-Galois}. 
For a finite group $G$, a commutative $R$-algebra $S$ with a $G$-action is Galois if the canonical maps 
\begin{equation*}
    R \;\too\; S^{hG}, \qquad S \otimes_R S \;\too\; S^G, \ x \otimes y \;\mapsto\; (g.x \cdot y)_{g \in G}
\end{equation*}
are isomorphisms.
We write $(-)^{hG}$ for fixed points, reserving $(-)^G$ for the $G$-indexed product.

These extensions correspond to $G$-torsors over $\Spec R$. When $R = K$ is a field, connected torsors recover classical Galois extensions. Disconnected torsors are allowed as well, corresponding to split Galois extensions, e.g. $S = K^G$.
In this case, the extensions are classified by continuous homomorphisms $\G_K \to G$, up to conjugation, without the surjectivity condition.

Deligne and Milne also classified these finite Galois extensions by exact symmetric monoidal functors~\cite[Theorem~3.2(b)]{Deligne-Milne-1982-Tannakian}.
A presentable version of this correspondence classifies Galois extensions of $R$ with group $G$ by colimit-preserving symmetric monoidal $R$-linear functors
\begin{equation*}
    \Mod_R(\Ab)^{\BG} \;\too\; \Mod_R(\Ab).
\end{equation*}

It is natural to ask how far this theory extends to general symmetric monoidal categories\footnote{Throughout this article, we will say `category' for `$(\infty,1)$-category'.}.
Rognes extended the definition of a Galois extension to general symmetric monoidal categories and $\EE_1$-groups in~\cite{Rognes-2008-Galois}, using the same two Galois conditions as above.
A $G$-action is equivalently a $\BG$-local system, and both Galois conditions can be expressed entirely in terms of the anima $\BG$.
The definition can therefore be extended to arbitrary animae, as was done in~\cite{BCSY-Fourier}.

Let $\cC$ be a presentably symmetric monoidal category, $R \in \CAlg(\cC)$ and $A$ an anima.
Write $\Delta \colon A \to A \times A$ for the diagonal and $A \colon A \to \pt$ for the terminal map.
We denote restriction along maps of animae by $(-)^*$ and left and right Kan extensions by $(-)_!$ and $(-)_*$ respectively.
We write 
\begin{equation*}
    \boxtimes_R \;\colon\; \Mod_R(\cC)^A \times \Mod_R(\cC)^A \;\too\; \Mod_R(\cC)^{A \times A}
\end{equation*}
for the external tensor product over $R$, given by $(M \boxtimes_R N)_{a,b} = M_a \otimes_R N_b$.

\begin{definition}[{\cite[Definition~2.25]{BCSY-Fourier}}]
\label{def:Galois}
    An $A$-local system $S$ of commutative $R$-algebras is an \emph{$A$-Galois extension of $R$} if:
    \begin{enumerate}[label=(G\arabic*), leftmargin=*]
        \item\label{cond:Galois-1}
        The mate $R \to A_* S$ of the unit $A^*R \to S$ is an isomorphism.
        \item\label{cond:Galois-2}
        The mate $S \boxtimes_R S \to \Delta_*S$ of the multiplication $\Delta^*(S \boxtimes_R S) \simeq S \otimes_R S \to S$ is an isomorphism.
    \end{enumerate}
    The extension is \emph{faithful} if $S \otimes_R - \colon \Mod_R(\cC) \to \Mod_R(\cC)^A$ is conservative.
    We denote the anima of faithful $A$-Galois extensions of $R$ by $\Gal_A(R; \cC)$.
\end{definition}

\begin{remark}
    When $A = \BG$ it recovers Rognes' definition of Galois extensions.
    When $A \simeq \bigsqcup_{a \in \pi_0 A} \BG_a$ is disconnected and $\cC$ is pointed, an $A$-Galois extension of $R$ consists of a decomposition
    \begin{equation*}
        R \;\simeq\; \bigsqcap_{a \in \pi_0 A} R_a
    \end{equation*}
    together with a $\BG_a$-Galois extension $S_a$ of $R_a$ for each $a \in \pi_0 A$.
    It is faithful when each $S_a$ is faithful over $R_a$ and base change induces an isomorphism
    \begin{equation*}
        \Mod_R(\cC) \;\isotoo\; \bigsqcap_{a \in \pi_0 A} \Mod_{R_a}(\cC).
    \end{equation*}
\end{remark}

\begin{remark}
    Classical Galois extensions are automatically faithful.
    Rognes introduced faithfulness as an additional condition for Galois extensions of ring spectra, leaving open whether it followed from the Galois conditions.
    Wieland subsequently gave a counterexample for commutative $\Fp$-ring spectra~\cite{Rognes-unfaithful}.
    Faithfulness plays a crucial role in the existing theory, particularly in Galois descent. We therefore restrict our attention to faithful Galois extensions. 
    The classification we prove below gives a further justification for this restriction.
\end{remark}

Replacing $\cC$ with $\Mod_R(\cC)$ we can always assume $R = \ounit_{\cC}$. We write $\Gal_{A}(\cC) \coloneq \Gal_{A}(\ounit_{\cC}; \cC)$.

Although the definition makes sense for arbitrary groups, the theory requires a suitable finiteness condition. Following~\cite[\S~3.4]{Rognes-2008-Galois}, we use dualizability.
An $\EE_1$-group $G$ is $\cC$-dualizable if the group algebra $\ounit_{\cC}[G]$ is dualizable as an object of $\cC$.
For an arbitrary anima $A$, the condition on the group becomes a condition on its path animae.
We say $A$ is \CO-dualizable if $\ounit_{\cC}[\Path_{ab}A]$ is dualizable in $\cC$ for every $a,b \in A$.

\begin{definition}
    Let $A$ be \CO-dualizable.
    For each $a \in A$, the local system $\triv_{A,a} \in \CAlg(\cC)^A$ given by
    \begin{equation*}
        (\triv_{A,a})_b \;\coloneq\; \ounit_{\cC}^{\Path_{ab}A}
    \end{equation*}
    is a faithful $A$-Galois extension, called the \emph{trivial $A$-Galois extension at $a$}.
    As $a$ varies, these assemble into an $(A \times A)$-local system.
    Viewed as an $A$-local system in $\cC^A$, it is a faithful $A$-Galois extension $\triv_A \in \Gal_A(\cC^A)$.
    
    When $A = \BG$, $\triv_{A,a}$ is the split algebra $\ounit_{\cC}^G$ with its left translation action, while $\triv_A$ carries the commuting left and right translation actions.
\end{definition}

For example, if $\cC$ is semiadditive (e.g.\ abelian or stable), finite groups are $\cC$-dualizable, and hence finite groupoids are \CO-dualizable.

The classification of finite Galois extensions of a field $K$ by $\B\G_K$ extends to a general commutative ring $R$ using its \'etale fundamental groupoid $\pi^{\et}_{\le 1} (\Spec R)$.
Mathew proved an analogous result for every presentably symmetric monoidal stable category $\cC$, classifying faithful finite Galois extensions by the profinite groupoid $\pi_{\le 1}^{\mrm{weak}}\cC$~\cite[Definition~6.8 and Proposition~6.13]{Mathew-2016-Galois}.

\subsection{Main results}

    Our main theorem extends the Deligne--Milne correspondence to general presentably symmetric monoidal categories. We use this result to construct Galois pro-animae generalizing the classical profinite absolute Galois groups and Mathew's weak Galois profinite groupoid.

    We fix a base category $\basecat \in \CAlg(\Pr)$ and let $\ModB \coloneq \Mod_{\basecat}(\Pr)$ be the category of $\basecat$-linear presentable categories.

    \begin{alphThm}[Corepresentability of Galois extensions~\cref{thm:Galois-corepresented}]
    \label{alphthm:Galois-co-representable}
        Let $\cC \in \CAlgB$ and $A$ be a \COB-dualizable anima.
        There is a natural isomorphism
        \begin{equation*}
            \begin{split}
                \MapBo(\basecat^A, \cC) & \;\isotoo\; \Gal_A(\cC), \\
                F & \;\mapsto\; F^A(\triv_A)
            \end{split}
        \end{equation*}
        where the left-hand side is the anima of colimit-preserving $\basecat$-linear symmetric monoidal functors, and $F^A \colon (\basecat^A)^A \to \cC^A$ applies $F$ pointwise.
    \end{alphThm}
    The inverse sends $R \in \Gal_A(\cC)$ to the functor
    \begin{equation*}
        F_R \;\colon\; \basecat^A \;\xtoo{\ounit_{\cC}[-]}\; \cC^A \;\xtoo{R \otimes -}\; \cC^A \;\xtoo{A_*}\; \cC.
    \end{equation*}

    The Galois conditions involve both tensor products and limits, which need not commute. Their compatibility is therefore an issue throughout the theory.
    For example, preservation of Galois extensions by colimit-preserving symmetric monoidal functors has been established in many cases, but the proofs require involved arguments~\cite[Lemma~7.1.1]{Rognes-2008-Galois},~\cite[Propositions~6.6 and~6.13]{Mathew-2016-Galois},~\cite[Proposition~3.12]{BMCSY-cycloshift}.

    Dualizability of the path animae of $A$ implies that every $A$-Galois extension is itself dualizable~(\cite[Proposition~6.2.1]{Rognes-2008-Galois},~\cref{prop:Gal-dualizable}).
    Its dual is an $A$-local system of commutative coalgebras.
    We call an $A$-local system of commutative coalgebras an \emph{$A$-coGalois extension} if it satisfies the conditions dual to \labelcref{cond:Galois-1,cond:Galois-2}~(\cref{def:coGal}).
    These conditions involve only tensor products and    colimits, which commute in a presentably symmetric monoidal category.

    Duality exchanges the second Galois and coGalois conditions.
    Extending results of Rognes, we show that the dual of an $A$-Galois extension satisfies the first coGalois condition precisely when the extension is faithful~(\cref{prop:Gal-coGal-iso}).
    Duality thus identifies faithful $A$-Galois extensions with $A$-coGalois extensions.
    To prove~\cref{alphthm:Galois-co-representable}, it therefore suffices to classify $A$-coGalois extensions. Writing $\coGal_A(\cC)$ for their anima, we prove the following theorem without any dualizability assumption on $A$.

    \begin{alphThmPrime}[\cref{prop:coGal-FunL-iso}]
    \label{thm:intro-coGal}
        Let $\cC \in \CAlgB$ and $A$ an anima.
        There is a natural isomorphism
        \begin{equation*}
            \begin{split}
                \MapBo(\basecat^A, \cC) & \;\isotoo\; \coGal_A(\cC), \\
                F & \;\mapsto\; F^A(\Delta_!A^*\ounit_{\basecat}).
            \end{split}
        \end{equation*}
    \end{alphThmPrime} 
    
    We prove~\cref{thm:intro-coGal} by first classifying colimit-preserving $\basecat$-linear oplax symmetric monoidal functors $\basecat^A \to \cC$.
    Since $\basecat^A \simeq \basecat \otimes \An^A$, it suffices to consider $\basecat = \An$.

    Since $\An^A$ is cartesian monoidal, the dual of~\cite[Theorem~2.4.3.18]{Lurie-HA} identifies oplax symmetric monoidal functors $\An^A \to \cC$ with functors $\An^A \to \coCAlg(\cC)$.
    Restricting to colimit-preserving functors, we obtain~(\cref{prop:oplax-monoidal-and-coalgebras})
    \begin{equation*}
        \FunLoplax(\An^A, \cC) \;\simeq\; \coCAlg(\cC)^A.
    \end{equation*}
    
    An oplax symmetric monoidal functor is symmetric monoidal if and only if its counit and comultiplication maps are isomorphisms.
    Under the above classification, these conditions correspond exactly to the coGalois conditions of the classifying local system of coalgebras.

    \cref{alphthm:Galois-co-representable} proves to be very useful in the development of Galois theory, allowing us to establish functoriality, descent, and relative Galois results with little additional work.
    As a first consequence, preservation of faithful Galois extensions under colimit-preserving symmetric monoidal functors is immediate.
    
    \begin{corollary}[\cref{cor:Galois-preserved-under-symmetric-monoidal}]
        Let $F \colon \cC \to \cD$ be a map in $\CAlgB$ and $A$ a \COB-dualizable anima. Then $F$ preserves faithful $A$-Galois extensions.
    \end{corollary}

    \subsubsection*{Galois descent}

    Classical Galois descent recovers modules over the base ring $R$ from modules over a $\BG$-Galois extension $R \to S$ equipped with a semilinear $G$-action.
    Equivalently, $\Mod_R(\Ab) \simeq \Mod_S(\Ab)^{hG}$.
    This result was extended to further settings~\cite[Theorem~9.4]{Mathew-2016-Galois},~\cite[Proposition~3.11]{BMCSY-cycloshift}.
    We use corepresentability to prove descent under the general assumption of \CO-dualizability.

    For an $A$-local system of commutative algebras $R$, write $\Mod_R(\cC)$ for the $A$-local system whose value at $a \in A$ is $\Mod_{R_a}(\cC)$.
    
    \begin{alphThm}[Galois descent~\cref{thm:Galois-descent}]
    \label{thm:intro-Galois-descent}
        Let $\cC \in \CAlg(\Pr)$ and $A$ be a \CO-dualizable anima.
        For any $R \in \Gal_A(\cC)$, the free-module functors induce a natural isomorphism
        \begin{equation*}
            \cC \;\isotoo\; A_* \Mod_R(\cC).
        \end{equation*}
    \end{alphThm}

    By~\cref{alphthm:Galois-co-representable}, the extension $R$ corresponds to a $\cC$-linear symmetric monoidal functor $F_R \colon \cC^A \to \cC$.
    The functors $a_! \colon \cC \to \cC^A$ exhibit $\cC^A$ as the colimit of an $A$-indexed diagram of copies of $\cC$, each equipped with the $\cC^A$-module structure given by $a^*$.
    Base change along $F_R$ turns it into the local system $\Mod_R(\cC)$~(\cref{prop:categorical-descent-pointwise}).
    Since base change preserves colimits, we obtain
    \begin{equation*}
        A_!\Mod_R(\cC) \;\simeq\; \cC \otimes_{\cC^A} \cC^A \;\simeq\; \cC.
    \end{equation*}
    Colimits in $\ModC$ along $\cC$-linear left adjoints are computed as limits in $\Cat$ of the corresponding right adjoints. 
    The functors $\Mod_{R_a}(\cC) \to \cC$ in this colimit presentation have right adjoints $R_a \otimes -$, giving the descent isomorphism.

    \subsubsection*{Relative Galois theory}

    Another aspect of Galois theory we wish to generalize is the description of intermediate extensions through subgroups and quotients of the Galois group.
    Let $E/K$ be a finite Galois extension of fields with group $G$.
    For a subgroup $H \subseteq G$, the extension $E/E^{hH}$ is Galois with group $H$.
    Writing $g \colon \BH \to \BG$ for the induced map, this is the restriction $g^*E$, regarded as an extension of $E^{hH} \simeq (\BH)_*g^*E$.

    If $H$ is normal, the extension $E^{hH}/K$ is Galois with group $G/H$.
    Writing $f \colon \BG \to \B(G/H)$ for the quotient map, the fixed field with its residual action is the right Kan extension $f_*E$.

    Using corepresentability, we prove analogous statements for general categories and \CO-dualizable animae, with $\cC$-affine maps replacing subgroup inclusions~(\cref{subsec:affineness}).
    For $\cC$ the category of modules over a nonzero commutative ring and a homomorphism $H \to G$ of finite groups, the induced map $\BH \to \BG$ is $\cC$-affine precisely when the homomorphism is injective.

    \begin{proposition}[\cref{thm:relative-quotient,thm:relative-subgroup}]
    \label{prop:intro-relative-Galois}
        Let $\cC$ be a presentably symmetric monoidal category, let $A,B,B'$ be \CO-dualizable animae and let $R \in \Gal_A(\cC)$.
        \begin{enumerate}
            \item Let $g \colon B' \to A$ be a $\cC$-affine map.
            Then $g^*R \in \Gal_{B'}(B'_*g^*R;\cC)$.
            \item Let $f \colon A \to B$ be any map.
            Then $f_*R \in \Gal_B(\cC)$.
        \end{enumerate}
    \end{proposition}
    
    Indeed, if $F_R \colon \cC^A \to \cC$ classifies $R$, then $F_R \circ f^*$ classifies $f_*R$.
    By affineness, $g_*$ identifies $\cC^{B'}$ with $\Mod_{g_*B'^*\ounit_{\cC}}(\cC^A)$.
    As $F_R(g_*B'^*\ounit_{\cC}) \simeq B'_*g^*R$, it induces a symmetric monoidal functor
    \begin{equation*}
        \cC^{B'} \;\too\; \Mod_{B'_*g^*R}(\cC),
    \end{equation*}
    which classifies $g^*R$.

    \subsubsection*{Galois animae}

    We use~\cref{alphthm:Galois-co-representable} to construct Galois pro-animae classifying Galois extensions.
    For this, we impose an additional affineness assumption.
    An anima $A$ is $\ModB$-affine~(\cite[Definition~2.1]{BCSY-Fourier}, see also~\cref{subsec:affineness}) if the limit functor identifies
    \begin{equation*}
        \ModB^A \;\isotoo\; \ModC[\basecat^A].
    \end{equation*}
    For example, if $\basecat$ is semiadditive, all $1$-trunacted $\pi$-finite animae (a.k.a.\ finite groupoids) are $\ModB$-affine~(\cite[Example~2.35]{BCSY-Fourier},~\cref{prop:affineness-and-categorification}).

    We fix a small full subcategory $\cT \subseteq \An$ consisting of truncated, \COB-dualizable and $\ModB$-affine animae, which is closed under finite limits.

    \begin{alphThm}[\cref{thm:pro-affine-anima}]
    \label{alphthm:prorepresentability}
        Let $\cC \in \CAlgB$. There exists a Galois pro-anima $\piT \cC \in \Pro(\cT)$ classifying $\cT$-Galois extensions: there is a natural isomorphism in $A \in \cT$
        \begin{equation*}
            \Gal_A(\cC) \;\simeq\; \Map_{\ProAn}(\piT \cC, A).
        \end{equation*}
    \end{alphThm}

    Prorepresentability of $\Gal_{(-)}(\cC) \colon \cT \to \An$ is equivalent to preservation of finite limits.
    By~\cref{alphthm:Galois-co-representable}, this functor agrees with $\MapBo(\basecat^{(-)},\cC)$.
    We therefore reduce to showing that $\basecat^{(-)}$ sends pullbacks in $\cT$ to pushouts in $\CAlgB$, which follows from~\cite[Theorem~2.38]{BCSY-Fourier}. 

    For $\basecat=\Ab$ and $\cT$ the category of finite groupoids, the category of $K$-vector spaces has Galois pro-anima $\piT\Vect_K \simeq \B\G_K$.
    More generally, for a commutative ring $R$, we recover the \'etale fundamental profinite groupoid $\piT \Mod_R(\Ab) \simeq \pi^{\et}_{\le 1}(\Spec R)$.

    For $\basecat = \Sp$ and $\cT$ the category of finite groupoids, our Galois pro-anima agrees with Mathew's weak Galois profinite groupoid, since both satisfy the same universal property~\cite[Definition~6.8 and Proposition~6.13]{Mathew-2016-Galois}.

    Categorification preserves the existing Galois theory.
    For a presentably symmetric monoidal $2$-category $\cD$ with $\End(\ounit_{\cD}) \simeq \cC$, we prove~(\cref{prop:piT-and-categorification}) that $\piT(\cC) \simeq \piT(\cD)$.
    Equivalently, $\Gal_A(\cC) \simeq \Gal_A(\cD)$ for every $A \in \cT$.
    For $\cD=\ModC$, this isomorphism sends $R$ to $\Mod_R(\cC)$.

    Categorification also allows us to enlarge the class of indexing animae.
    For $\cC=\Mod_R(\Ab)$, passing from finite groupoids to $2$-truncated $\pi$-finite animae gives $\cD$ a $2$-truncated pro-$\pi$-finite Galois anima whose $1$-truncation is $\pi^{\et}_{\le 1}(\Spec R)$.
    For example, the symmetric monoidal category $\sVect_{\CC}$ of complex super vector spaces is a nontrivial $\B^2\ZZ/2$-Galois extension of $\Vect_{\CC}$ in $\ModC[\Vect_{\CC}]$~(\cite[Theorem~2.11]{Johson-Freyd-2017-sVect},~\cite[Proposition~5.1.4 and Example~1.3.9]{Keidar-Ragimov-2025-twisted}).

\subsection{Organization}
    In~\cref{sec:corepresentability}, we prove~\cref{alphthm:Galois-co-representable}.
    We classify oplax symmetric monoidal functors by local systems of coalgebras in~\cref{subsec:oplax-functors}, identify symmetric monoidal functors with coGalois extensions in~\cref{subsec:cogal}, and use duality to recover faithful Galois extensions in~\cref{subsec:norm-and-coGalois}.

    In~\cref{sec:affine-galois-theory}, we use corepresentability to prove Galois descent and relative Galois results.
    We establish the required preliminaries in~\cref{subsec:trivial-extensions,subsec:affineness}, then prove descent in~\cref{subsec:Galois-descent} and the relative Galois results in~\cref{subsec:relative-Galois}.

    In~\cref{sec:Galois-anima}, we develop the theory of Galois pro-animae.
    We prove~\cref{alphthm:prorepresentability} in~\cref{subsec:prorepresentability} and study their behavior under categorification in~\cref{subsec:categorification}.
    In~\cref{subsec:pro-Galois}, we extend both classification theorems to pro-Galois extensions.

\subsection{Acknowledgments}
    I thank Shay Ben-Moshe, Itamar Mor, Niko Naumann, Leor Neuhauser, Shaul Ragimov, Tomer Schlank, and Lior Yanovski for helpful discussions. 
    I am especially grateful to Shaul Ragimov for many hours of conversation about this and related projects.
    I thank Shay Ben-Moshe, Shaul Ragimov and Lior Yanovski for helpful comments on previous versions. 
    
    This work was partially supported by the CRC 1785 Generalised Motivic Methods in Geometry (Universit\"at Regensburg, funded by the DFG).

    \subsubsection*{AI disclosure.}  
    All ideas, results and proofs were obtained without the use of AI.
    We used ChatGPT for help with writing, locating references, and proofreading. It identified some minor errors and suggested a simplification of the proof of~\cref{prop:oplax-monoidal-and-coalgebras}.


\section{Galois Corepresentability}
\label{sec:corepresentability}
In this section, we prove the corepresentability of Galois extensions~\cref{alphthm:Galois-co-representable}.

The proof is organized into three steps. 
In~\cref{subsec:oplax-functors}, we classify colimit-preserving oplax symmetric monoidal functors $\An^A \to \cC$ by $A$-local systems of coalgebras in $\cC$.
In~\cref{subsec:cogal}, we show that strong monoidality corresponds to \quotes{coGalois} conditions, dual to those of~\cref{def:Galois}.
Finally, in~\cref{subsec:norm-and-coGalois}, under a dualizability assumption on $A$, we show that duality identifies coGalois extensions with faithful Galois extensions.

We fix a base category $\basecat \in \CAlg(\Pr)$ and work in the category of presentable $\basecat$-linear categories $\ModB = \Mod_{\basecat}(\Pr)$.
We also fix $\cC\in\CAlgB$.

Let $f \colon A \to B$ be a map of animae. We denote by $f^*$ the $\basecat$-linear, symmetric monoidal functor
\begin{equation*}
    f^* \;\colon\; \cC^B \;\too\; \cC^A,
\end{equation*}
leaving the category $\cC$ implicit.
It admits both a left and a right adjoint, which we denote by $f_!$ and $f_*$, respectively.

We fix notation for several maps used throughout the paper:
\begin{notation}
    Let $f \colon A \to B$ and $f' \colon A' \to B$ be maps of animae.
    \begin{itemize}
        \item The terminal map $A \to \pt$ will also be denoted by $A$, using the identification $\An \simeq \An_{/\pt}$.
        \item The projection maps of the fibered product $A \times_B A'$ are denoted by $\pr_f$ and $\pr_{f'}$.
        \item The projection maps of the product $A \times_B A$ are called $\pr_1, \pr_2$.
        \item The diagonal map is denoted $\Delta_f \colon A \to A \times_B A$. If $f$ is clear from context, it will just be called~$\Delta$.
    \end{itemize}
\end{notation}

\begin{notation}
    Let $\cC \in \CAlgB$ be a $\basecat$-algebra. We denote by
    \begin{equation*}
        \ounit_{\cC}[-] \;\colon\; \basecat \;\too\; \cC
    \end{equation*}
    the unique $\basecat$-linear symmetric monoidal functor.
    For $A \in \basecat$ and $X \in \cC$, we write $X[A] \coloneq X \otimes \ounit_{\cC}[A]$.

     For $\basecat = \An$, $X[A]$ is the constant $A$-shaped colimit of $X$.
\end{notation}

\subsection{Classification of oplax symmetric monoidal functors}
\label{subsec:oplax-functors}
    We first recall the classification of colimit-preserving $\basecat$-linear functors $\basecat^A \to \cC$.
    The category $\basecat^A$ is freely generated, as a presentable $\basecat$-linear category, by the following path local system.

    \begin{definition}
        Consider the path local system $\Path A \in \An^{A \times A}$ defined by
        \begin{equation*}
            \Path A \;\coloneq\; \Delta_! A^* \pt.
        \end{equation*}
        Equivalently, at $(a,b)$ it is the anima $\Path_{ab}A$ of paths between $a$ and $b$ in $A$.
    \end{definition}

    The local system $\ounit_{\basecat}[\Path A] \simeq \Delta_! A^* \ounit_{\basecat} \in (\basecat^A)^A$ is the universal $A$-local system in the following sense:
    \begin{proposition}[{\cite[Corollary~2.30]{Cnossen-2023-twisted-ambi}, see also~\cite[Theorem~5.1.5.6]{Lurie-HTT}}]
    \label{prop:classification-of-functors}
        Let $\cC\in\CAlgB$ and $A\in\An$. The evaluation at $\ounit_{\basecat}[\Path A]$ functor
        \begin{equation*}
            \FunB(\basecat^A, \cC) \;\xtoo{(-)^A}\; \FunB(\basecat^{A\times A}, \cC^A) \;\xtoo{\ev_{\ounit_{\basecat}[\Path A]}}\; \cC^A
        \end{equation*}
        is an isomorphism of categories.
        The inverse sends a local system $X\in \cC^A$ to the functor
        \begin{equation*}
            F_X \;\colon\; \basecat^A \;\xtoo{\ounit_{\cC}[-]}\; \cC^A \;\xtoo{X\otimes -}\; \cC^A \;\xtoo{A_!}\; \cC.
        \end{equation*}
    \end{proposition}
    
    As $\basecat^A \simeq \An^A \otimes \basecat$, the $\basecat$-linearity identifies $\basecat$-linear colimit-preserving oplax symmetric monoidal functors $\basecat^A \to \cC$ with colimit-preserving oplax symmetric monoidal functors $\An^A \to \cC$. This case is simpler as $\An^A$ is equipped with the cartesian monoidal structure. Therefore, by the dual statement to~\cite[Theorem~2.4.3.18]{Lurie-HA}, we have
    \begin{equation*}
        \Funooplax(\An^A, \cC) \;\simeq\; \Fun(\An^A, \coCAlg(\cC)).
    \end{equation*}

    We define oplax symmetric monoidal colimit-preserving functors as oplax symmetric monoidal functors whose underlying functors are colimit-preserving. Formally:
    \begin{definition}
        Let $\cC \in \CAlg(\Pr)$. Consider the functors
        \begin{equation*}
            \Fun(\An^{(-)}, \cC),\ \FunL(\An^{(-)}, \cC),\ \Funooplax(\An^{(-)}, \cC) \;\colon\; \An \;\too\; \Cat.
        \end{equation*}
        The functoriality for a map $f \colon A \to B$ is by precomposition with $f^* \colon \An^B \to \An^A$.

        Define $\FunLoplax(\An^{(-)}, \cC) \colon \An \to \Cat$ as the pullback 
        \begin{equation*}
            \FunLoplax(\An^{(-)}, \cC) \;\coloneq\; \FunL(\An^{(-)}, \cC) \times_{\Fun(\An^{(-)}, \cC)} \Funooplax(\An^{(-)}, \cC).
        \end{equation*}
    \end{definition}
    
    \begin{proposition}
    \label{prop:oplax-monoidal-and-coalgebras}
        The isomorphism of~\cref{prop:classification-of-functors} for $\basecat=\An$ extends to an isomorphism of categories
        \begin{equation*}
            \FunLoplax(\An^A, \cC) \;\isotoo\; \coCAlg(\cC)^A.
        \end{equation*}
    \end{proposition}
    \begin{proof}
        By definition, $\FunLoplax(\An^A, \cC) = \FunL(\An^A, \cC) \times_{\Fun(\An^A, \cC)} \Funooplax(\An^A, \cC)$. The last term is identified with $\Fun(\An^A, \coCAlg(\cC))$, and the forgetful map to $\Fun(\An^A, \cC)$ is obtained by composing with $\coCAlg(\cC) \to \cC$.

        Therefore, a colimit-preserving oplax symmetric monoidal functor $\An^A \to \cC$ is equivalent to a functor $\An^A \to \coCAlg(\cC)$ such that the composition $\An^A \to \coCAlg(\cC) \to \cC$ commutes with colimits. 
        Since $\coCAlg(\cC) \to \cC$ preserves and detects colimits by the dual of~\cite[Corollary~3.2.2.5]{Lurie-HA}, this is equivalent to $\An^A \to \coCAlg(\cC)$ preserving colimits.
        Therefore
        \begin{equation*}
            \FunLoplax(\An^A, \cC) \;\simeq\; \FunL(\An^A, \coCAlg(\cC)).
        \end{equation*}
        The evaluation at $\Path A$ functor identifies the latter with $\coCAlg(\cC)^A$. 
    \end{proof}

\subsection{Classification of symmetric monoidal functors}
\label{subsec:cogal}
    In the previous subsection, we classified colimit-preserving oplax symmetric monoidal functors $\An^A \to \cC$ as corresponding to $A$-local systems of coalgebras.
    An oplax symmetric monoidal functor $F \colon \An^A \to \cC$ is symmetric monoidal if it satisfies
    \begin{enumerate}[label=(SM\arabic*), leftmargin=*]
        \item \label{cond:SM-1}
        The counit map $F(A^* \pt) \to \ounit_{\cC}$ is an isomorphism.
        \item \label{cond:SM-2}
        The comultiplication natural transformation $F(X \times Y) \to FX \otimes FY$ is a natural isomorphism.
    \end{enumerate}
    In this subsection, we classify the coalgebras that correspond to symmetric monoidal functors. These coalgebras satisfy conditions dual to~\cref{def:Galois}.
        
    \begin{definition}
    \label{def:coGal}
        Let $\cC\in \CAlg(\Pr)$ and let $A$ be an anima. A local system of cocommutative coalgebras $Q \in \coCAlg(\cC)^A$ is called an $A$-coGalois extension (of the unit object $\ounit_{\cC}$) if it satisfies the following two properties:
        \begin{enumerate}[label=(cG\arabic*), leftmargin=*]
            \item \label{cond:coGal-1}
            The mate $A_! Q \to \ounit_{\cC}$ of the counit $Q \to A^* \ounit_{\cC}$ is an isomorphism in $\coCAlg(\cC)$.
            \item \label{cond:coGal-2}
            The mate $\Delta_! Q \to Q \boxtimes Q$ of the comultiplication $Q \to Q \otimes Q \simeq \Delta^*(Q \boxtimes Q)$ is an isomorphism in $\coCAlg(\cC)^{A\times A}$.
        \end{enumerate}
        We will denote the anima of $A$-coGalois extensions by $\coGal_A(\cC)$.
    \end{definition}
        
    \begin{example}
    \label{exm:coGal-in-spc}
        Let $B$ be an anima, then 
        \begin{equation*}
            \Map(B, A) \;\isotoo\; \coCAlg_A(\An^B), \qquad f \;\longmapsto\; (f,\id_B)_!B^*\pt.
        \end{equation*}
    \end{example}
    \begin{proof}
        Write $\cC = \An^B$. We have an isomorphism $\coCAlg(\cC) \simeq \cC$, where the coalgebra structure of an object is given by its diagonal map. By unstraightening~\cite[Theorem~3.2.0.1]{Lurie-HTT},
        \begin{equation*}
            \cC \;\simeq\; \overCat{\An}{B}
            \qand
            \cC^A \;\simeq\; \overCat{\An}{A\times B}.
        \end{equation*}

        Under these identifications, the functors
        \begin{equation*}
            A_! \;\colon\; \cC^A \;\too\; \cC,
            \qquad
            \Delta_! \;\colon\; \cC^A \;\too\; \cC^{A\times A}
        \end{equation*}
        are given by post-composition with
        \begin{equation*}
            \pr_B \;\colon\; A\times B \;\too\; B,
            \qquad
            \Delta_A\times\id_B \;\colon\; A\times B
            \;\too\; A\times A\times B,
        \end{equation*}
        respectively.
        
        Let $q \colon X\to A\times B$ correspond to an object
        $Q\in\coCAlg(\cC)^A$. The object $A_!Q$ corresponds to the composite
        \begin{equation*}
            X \;\xtoo{q}\; A\times B \;\xtoo{\pr_B}\; B.
        \end{equation*}
        Therefore, $Q$ satisfies condition~\labelcref{cond:coGal-1} if and only if this composite is the terminal object of $\overCat{\An}{B}$, or equivalently, if $\pr_B\circ q\colon X\to B$ is an isomorphism. 
    
        On the other hand, the map in condition~\labelcref{cond:coGal-2} corresponds to the relative diagonal
        \begin{equation*}
            X \;\too\; X\times_B X.
        \end{equation*}
        Hence condition~\labelcref{cond:coGal-2} holds automatically whenever condition~\labelcref{cond:coGal-1} holds.
    
        After identifying $X$ with $B$, the map $q$ is the graph
        \begin{equation*}
            (f,\id_B) \;\colon\; B \;\too\; A\times B
        \end{equation*}
        of a map $f\colon B\to A$. We therefore obtain the required isomorphism
        \begin{equation*}
            \coGal_A(\An^B) \;\simeq\; \Map(B,A). \qedhere
        \end{equation*}
    \end{proof}
    
    \begin{example}
    \label{exm:universal-coGalois}
        As a special case of the previous example, the path local system $\Path A = \Delta_! A^* \pt$ is an $A$-coGalois extension in $\An^A$ corresponding to the identity map of $A$. 
        \cref{prop:coGal-FunL-iso} shows that it is the universal $A$-coGalois extension.
    \end{example}

    We will use the following compatibility with pullback and left Kan extension.

    \begin{lemma}
    \label{lem:lower-shriek-functoriality}
        Let $F \colon \cC \to \cD$ be a morphism in $\Pr$, and let $f \colon A \to B$ be a map of animae. Then the following diagrams commute:
        \begin{equation*}
            \begin{tikzcd}
                {\cC^A} && {\cD^A} \\
                \\
                {\cC^B} && {\cD^B,}
                \arrow["F^A", from=1-1, to=1-3]
                \arrow["{f_!}", from=1-1, to=3-1]
                \arrow["{f_!}", from=1-3, to=3-3]
                \arrow["F^B", from=3-1, to=3-3]
            \end{tikzcd}
            \qand
            \begin{tikzcd}
                {\cC^B} && {\cD^B} \\
                \\
                {\cC^A} && {\cD^A.}
                \arrow["F^B", from=1-1, to=1-3]
                \arrow["{f^*}", from=1-1, to=3-1]
                \arrow["{f^*}", from=1-3, to=3-3]
                \arrow["F^A", from=3-1, to=3-3]
            \end{tikzcd}
        \end{equation*}
    \end{lemma}
    \begin{proof}
        $F$ defines a functor $\pt \to \Pr^{[1]}$. Consider its left Kan extension
        \begin{equation*}
            \begin{tikzcd}
                \pt && {\Pr^{[1]}.} \\
                {\An}
                \arrow["F", from=1-1, to=1-3]
                \arrow[from=1-1, to=2-1]
                \arrow[dashed, from=2-1, to=1-3]
            \end{tikzcd}
        \end{equation*}
        This takes an anima $A$ to the map $F[A] \colon \cC[A] \to \cD[A]$, and a map of animae $f \colon A \to B$ to the required square with $(-)_!$-functoriality, using the identifications $\cC[A] \simeq \cC^A$, $\cD[A] \simeq \cD^A$.
        Similarly, consider its right Kan extension
        \begin{equation*}
            \begin{tikzcd}
                \pt && {\Pr^{[1]}.} \\
                {\An\op}
                \arrow["F", from=1-1, to=1-3]
                \arrow[from=1-1, to=2-1]
                \arrow[dashed, from=2-1, to=1-3]
            \end{tikzcd}
        \end{equation*}
        This takes an anima $A$ to the map $F^A \colon \cC^A \to \cD^A$, and a map of animae $f \colon A \to B$ to the required square with $(-)^*$-functoriality.
    \end{proof}

    \begin{theorem}
    \label{prop:coGal-FunL-iso}
        Let $F \colon \An^A \to \cC$ be a colimit-preserving oplax symmetric monoidal functor and let $Q \in \coCAlg(\cC)^A$ be its image under the isomorphism of~\cref{prop:oplax-monoidal-and-coalgebras}. Then
        \begin{enumerate}
            \item $F$ satisfies~\labelcref{cond:SM-1} if and only if $Q$ satisfies~\labelcref{cond:coGal-1}.
            \item $F$ satisfies~\labelcref{cond:SM-2} if and only if $Q$ satisfies~\labelcref{cond:coGal-2}.
        \end{enumerate}
        In particular, the isomorphism of~\cref{prop:oplax-monoidal-and-coalgebras} restricts to an isomorphism of animae
        \begin{equation*}
            \MapBo(\basecat^A, \cC) \;\simeq\; \MapCo[](\An^A, \cC) \;\isotoo\; \coGal_A(\cC).
        \end{equation*}
    \end{theorem}
    \begin{proof}
        Writing the isomorphism explicitly, we have
        \begin{equation*}
            Q = F^A(\Path A), \qquad F = A_!(Q \otimes \ounit_{\cC}[-]).
        \end{equation*}
        Using~\cref{lem:lower-shriek-functoriality} for the terminal map $A \to \pt$, together with the fact that $\Path A$ is $A$-coGalois in $\An^A$ (\cref{exm:universal-coGalois}), we identify the counit map of $F$ with the counit map of $Q$:
        \begin{equation*}
            A_! Q \;=\; A_! F^A(\Path A) \;\simeq\; F(A_! \Path A) \;\simeq\; F(A^* \pt) \;\too\; \ounit_{\cC}.
        \end{equation*}
        Therefore, $F$ satisfies~\labelcref{cond:SM-1} if and only if $Q$ satisfies~\labelcref{cond:coGal-1}.

        Write $\eta_{X,Y} \colon F(X \times Y) \to F(X) \otimes F(Y)$ for the comultiplication of $F$. Both sides preserve colimits in each variable. Therefore, by~\cref{prop:classification-of-functors}, condition~\labelcref{cond:SM-2} holds if and only if the map
        \begin{equation} \label{eq:oplax-monoidal-map} \tag{$\star$}
            \eta_{\Path A, \Path A} \;\colon\; F^{A \times A}(\Path A \boxtimes^{\An^A} \Path A) \;\too\; F^A(\Path A) \boxtimes^{\cC} F^A(\Path A)
        \end{equation}
        is an isomorphism. Here
        \begin{equation*}
            \boxtimes^{\An^A} \;\colon\; \An^{A \times A} \times \An^{A \times A} \;\too\; \An^{A \times A \times A}
        \end{equation*}
        is the external product of $(\An^A)^A$, given by $(X\boxtimes^{\An^A}Y)(a,b,x)=X(a,x)\times Y(b,x)$, where $x$ is the inner coordinate. Similarly,
        \begin{equation*}
            \boxtimes^{\cC} \;\colon\; \cC^A \times \cC^A \;\too\; \cC^{A \times A}
        \end{equation*}
        is the external tensor product of $\cC^A$.

        As $\Path A$ is $A$-coGalois, 
        \begin{equation*}
            \Path A \boxtimes^{\An^A} \Path A \;\simeq\; \Delta_! \Path A \qin \An^{A \times A \times A}.
        \end{equation*}
        Therefore, by~\cref{lem:lower-shriek-functoriality},
        \begin{equation*}
            F^{A \times A} (\Path A \boxtimes^{\An^A} \Path A) \;\simeq\; \Delta_! F^A(\Path A)
        \end{equation*}
        and the map~\labelcref{eq:oplax-monoidal-map} is identified with the mate of the comultiplication of $Q$:
        \begin{equation*}
            \Delta_! Q \;\too\; Q \boxtimes^{\cC} Q.
        \end{equation*}
        Therefore, \labelcref{eq:oplax-monoidal-map} is an isomorphism if and only if $Q$ satisfies~\labelcref{cond:coGal-2}.
        The first isomorphism in the statement follows by extension of scalars, using $\basecat^A \simeq \basecat \otimes \An^A$ in $\CAlg(\Pr)$ (\cite[Corollary~3.4.1.7 and Proposition~3.2.4.7]{Lurie-HA}).
    \end{proof}

    This identification also gives an immediate proof of the following folklore result.
    \begin{corollary}
        The functor $\An^{(-)}\colon \An\op \to \CAlg(\Pr)$ is fully faithful.
    \end{corollary}
    \begin{proof}
        This follows from~\cref{prop:coGal-FunL-iso} and~\cref{exm:coGal-in-spc}.
    \end{proof}

\subsection{Dualizability and Galois extensions}
\label{subsec:norm-and-coGalois}
    \cref{prop:coGal-FunL-iso} identifies colimit-preserving $\basecat$-linear symmetric monoidal functors $\basecat^A \to \cC$ with $A$-coGalois extensions.
    Under a dualizability assumption, we identify the anima of $A$-coGalois extensions with the anima of faithful $A$-Galois extensions by sending a coGalois extension $Q$ to its dual $Q\dual$.

    \begin{definition}
        Let $A \in \An$ and $\basecat \in \CAlg(\Pr)$. We say that $A$ is \COB-dualizable if $\ounit_{\basecat}[\Path A] \in \basecat^{A \times A}$ is dualizable. Dualizability in $\basecat^{A\times A}$ is detected pointwise.
        In particular:
        \begin{enumerate}
            \item If $A = \BG$ is connected, then $A$ is \COB-dualizable if and only if $\ounit_{\basecat}[G]$ is dualizable.
            \item If $A = \emptyset$, then $A$ is \COB-dualizable as $\basecat^{A \times A} = \pt$.
            \item If $A$ is a nonempty non-connected anima, then it is \COB-dualizable if and only if $\emptyset$ is dualizable (equivalently, $\basecat$ is pointed) and $\ounit_{\basecat}[\Omega_a A]$ is dualizable for every $a \in A$.
        \end{enumerate}
    \end{definition}

    Note that if $A$ is \COB-dualizable and $\cC \in \CAlgB$, then $A$ is also \CO-dualizable.

    \begin{corollary}
    \label{cor:coGal-dualizable}
        Let $A$ be \COB-dualizable, $\cC \in \CAlgB$ and let $Q \in \coGal_A(\cC)$. Then $Q \in \cC^A$ is dualizable.
    \end{corollary}
    \begin{proof}
        Let $F_Q \colon \basecat^A \to \cC$ be the corresponding $\basecat$-linear symmetric monoidal functor. Then $Q$ is the image of $\ounit_{\basecat}[\Path A]$ under the symmetric monoidal functor
        \begin{equation*}
            \basecat^{A \times A} \;\simeq\; (\basecat^A)^A \;\xtoo{F_Q^A}\; \cC^A.
        \end{equation*}
        Since $A$ is \COB-dualizable, $\ounit_{\basecat}[\Path A] \in \basecat^{A \times A}$ is dualizable, and therefore so is $Q$ as its image under a symmetric monoidal functor.
    \end{proof}
    
    \subsubsection*{Projection formulas}
    
    We recall some known facts about projection formulas, which we employ to compute duals of Kan extensions. All of these can be found in~\cite{Fausk-Hu-May-2003-6FF}. 

    \begin{definition}
        Let $f^* \colon \cC \to \cD$ be a symmetric monoidal functor. Let $X \in \cC$ and $Y \in \cD$.
        \begin{enumerate}
            \item Assume $f^*$ admits a right adjoint $f_*$. Then $f_*$ is lax monoidal, and using the unit of the adjunction $u\colon \id \To f_* f^*$, we get the \emph{right projection map}
            \begin{equation*}
                \Rproj \;\colon\; X \otimes f_* Y \;\xtoo{u}\; f_* f^* X \otimes f_* Y \;\too\; f_*(f^*X \otimes Y).
            \end{equation*}

            \item Assume $f^*$ admits a left adjoint $f_!$. Then $f_!$ is oplax monoidal, and using the counit of the adjunction, $c \colon f_! f^* \To \id$, we get the \emph{left projection map}
            \begin{equation*}
                \Lproj \;\colon\; f_!(f^* X \otimes Y) \;\too\; f_! f^*X \otimes f_! Y \;\xtoo{c}\; X \otimes f_! Y.
            \end{equation*}
        \end{enumerate}
        We say that $f^*$ satisfies the right (resp.\ left) projection formula at $X \in \cC$ if, for any $Y \in \cD$, the right (resp.\ left) projection map is an isomorphism. We say that $f^*$ satisfies the right (resp.\ left) projection formula if it satisfies it at every $X\in \cC$. 
    \end{definition}

    \begin{lemma}[{\cite[Propositions~3.2 and~4.5]{Fausk-Hu-May-2003-6FF}}]
    \label{lem:projection-for-dbl}
        If $X \in \cC$ is dualizable, then $f^*$ satisfies the right and left projection formulas at $X$, whenever the corresponding adjoints exist.
    \end{lemma}
    \begin{proof}
        Let $Y \in \cD$. For any $Z \in \cC$,
        \begin{equation*}
            \begin{split}
                \Map_{\cC}(Z, X \otimes f_* Y) 
                & \;\simeq\; \Map_{\cC}(X\dual \otimes Z, f_* Y) \\
                & \;\simeq\; \Map_{\cD} (f^*(X\dual \otimes Z), Y) \\
                & \;\simeq\; \Map_{\cD}((f^* X)\dual \otimes f^* Z, Y) \\
                & \;\simeq\; \Map_{\cD}(f^* Z, f^* X \otimes Y) \\
                & \;\simeq\; \Map_{\cC}(Z, f_*(f^* X \otimes Y)).
            \end{split}
        \end{equation*}
        Tracing through these isomorphisms, the resulting map is induced by the right projection map. Therefore, by the Yoneda lemma, the right projection map is an isomorphism. The proof for the left projection map is similar.
    \end{proof}

    \begin{remark}
    \label{rmrk:left-projection-to-lax-closed}
        Assume that $\cC,\cD$ are closed monoidal and $f^* \colon \cC \to \cD$ admits a left adjoint $f_!$. The left projection gives $f^*$ a \emph{lax closed monoidal} structure, i.e.\ a natural transformation
        \begin{equation*}
            f^* \hom_{\cC}(-, -) \;\Too\; \hom_{\cD}(f^* -, f^* -) \qin \Fun(\cC\op \times \cC, \cD),
        \end{equation*}
        where $\hom_{\cC}$, $\hom_{\cD}$ are the internal hom functors of $\cC$, $\cD$, respectively.
        One can construct it on objects $X,Y \in \cC$ using the Yoneda lemma: for each $Z \in \cD$, we get a map
        \begin{equation*}
            \begin{split}
                \Map_{\cD}(Z, f^* \hom_{\cC}(X,Y)) 
                & \;\simeq\; \Map_{\cC}(f_! Z, \hom_{\cC}(X, Y)) 
                \;\simeq\; \Map_{\cC}(X \otimes f_! Z, Y) \\
                & \xtoo{- \circ \Lproj} \Map_{\cC}(f_!(f^* X \otimes Z), Y)
                \;\simeq\; \Map_{\cD}(f^*X \otimes Z, f^* Y) \\
                & \;\simeq\; \Map_{\cD}(Z, \hom_{\cD}(f^* X, f^* Y)).
            \end{split}
        \end{equation*}
        When this natural transformation is an isomorphism, we say that $f^*$ is closed monoidal.
    \end{remark}

    \begin{corollary}[{\cite[Proposition~4.3]{Fausk-Hu-May-2003-6FF}}]
    \label{lem:left-projection}
        Assume $\cC, \cD$ are closed monoidal and $f^* \colon \cC \to \cD$ admits a left adjoint $f_!$. Then $f^*$ is closed monoidal if and only if $f_!$ satisfies the left projection formula.
    \end{corollary}
    \begin{proof}
        The \quotes{if} direction follows from~\cref{rmrk:left-projection-to-lax-closed}.
        Assume $f^*$ is closed monoidal. Then, again by~\cref{rmrk:left-projection-to-lax-closed}, for any $X, Y\in \cC$ and $Z \in \cD$,
        \begin{equation*}
            \Map_{\cC}(X \otimes f_! Z, Y) \;\xtoo{- \circ \Lproj}\; \Map_{\cC}(f_!(f^* X \otimes Z), Y)
        \end{equation*}
        is an isomorphism. By the Yoneda lemma, the projection map is an isomorphism.
    \end{proof}

    \begin{example}
    \label{exm:f!-projection-formula}
        Let $f \colon A \to B$ be a map of animae and $\cC \in \CAlg(\Pr)$. Then $f_!$ satisfies the left projection formula, or equivalently, $f^*$ is closed monoidal. Indeed, we get a map for any $X, Y \in \cC^B$
        \begin{equation*}
            f^* \hom(X, Y) \;\too\; \hom(f^* X, f^* Y),
        \end{equation*}
        which after composing with $a^*$ is an isomorphism for any $a \in A$. As the $a^*$ are jointly conservative, this map is an isomorphism.
    \end{example}

    \begin{corollary}
    \label{cor:internal-adjoints-maps-of-animae}
        Let $f \colon A \to B$ be a map of animae and $\cC \in \CAlg(\Pr)$. Then for $X \in \cC^A$, $Y \in \cC^B$
        \begin{equation*}
            \begin{split}
                & \hom(f_! X, Y) \;\simeq\; f_*\hom(X, f^* Y) \qin \cC^B, \\
                & \hom(Y, f_* X) \;\simeq\; f_*\hom(f^* Y, X) \qin \cC^B.
            \end{split}
        \end{equation*}
        In particular, the weak dual of $f_! X$ is
        \begin{equation*}
            \hom(f_! X, B^* \ounit_{\cC}) \;\simeq\; f_* \hom(X, A^* \ounit_{\cC}).
        \end{equation*}
    \end{corollary}
    \begin{proof}
        Let $Z \in \cC^B$. Then, using the closed monoidality of $f^*$,
        \begin{equation*}
            \begin{split}
                \Map(Z, \hom(f_! X, Y)) 
                & \;\simeq\; \Map(f_! X, \hom(Z, Y)) \\
                & \;\simeq\; \Map(X, \hom(f^* Z, f^* Y)) \\
                & \;\simeq\; \Map(f^* Z, \hom(X, f^* Y)) \\
                & \;\simeq\; \Map(Z, f_* \hom(X, f^* Y)).
            \end{split}
        \end{equation*}
        Similarly,
        \begin{equation*}
            \begin{split}
                \Map(Z, f_* \hom(f^* Y, X))
                & \;\simeq\; \Map(f^* Y \otimes f^* Z, X) \\
                & \;\simeq\; \Map(Y \otimes Z, f_* X) \\
                & \;\simeq\; \Map(Z, \hom(Y, f_* X)).
            \end{split}
        \end{equation*}
    \end{proof}

    In particular, the weak dual of the universal coGalois extension is $\ounit_{\basecat}^{\Path A} \simeq \Delta_* A^* \ounit_{\basecat}$.

    \begin{lemma}
    \label{lem:Y-and-f!Y-dual-projection}
        Let $f \colon A \to B$ be a map of animae and $\cC \in \CAlg(\Pr)$. Let $Y \in \cC^A$ be such that $Y$ and $f_!(Y\dual)$ are dualizable. Then for any $X \in \cC^B$, the right projection map is an isomorphism
        \begin{equation*}
            \Rproj \;\colon\; X \otimes f_* Y \;\isotoo\; f_*(f^* X \otimes Y).
        \end{equation*}
    \end{lemma}
    \begin{proof}
        By~\cref{cor:internal-adjoints-maps-of-animae}, $f_*Y$ is dualizable and $f_!(Y\dual) \simeq (f_*Y)\dual$. Let $Z \in \cC^B$. Using the closed monoidality of $f^*$,
        \begin{equation*}
            \begin{split}
                \Map(Z, X \otimes f_* Y) 
                &\;\simeq\; \Map(Z \otimes f_!(Y\dual), X) \\
                & \;\simeq\; \Map(f_!(Y\dual), \hom(Z, X)) \\
                & \;\simeq\; \Map(Y\dual, \hom(f^* Z, f^* X)) \\
                & \;\simeq\; \Map(f^* Z, f^* X \otimes Y) \\
                & \;\simeq\; \Map(Z, f_*(f^* X \otimes Y)).
            \end{split}
        \end{equation*}
        Tracing through these isomorphisms, the resulting map is induced by the right projection map.
    \end{proof}

    \subsubsection*{The norm map and duality}

    Let $f \colon A \to B$ be a map of animae. The dualizing object of $f$ is defined as the local system
    \begin{equation*}
        \DD_f \;\coloneq\; (\pr_1)_* (\Delta_f)_! A^* \ounit_{\cC} \;\simeq\; (\pr_2)_* (\Delta_f)_! A^* \ounit_{\cC} \qin \cC^A.
    \end{equation*}
    In particular, $\DD_A = {(\pr_1)}_* \ounit_{\cC}[\Path A]$.
    There is always a natural transformation, called the norm, relating the colimit and limit functors along the fibers of $f$:
    \begin{equation*}
        \Nm_f \;\colon\; f_!( - \otimes \DD_f) \;\Too\; f_* \qin \Fun(\cC^A, \cC^B).
    \end{equation*}
    The mate of $\Nm_f$ is induced, via Beck--Chevalley and the left projection formula, by the counit $\pr_1^*\DD_f \to {(\Delta_f)}_!A^*\ounit_{\cC}$; see~\cite[Lemma~3.6]{Cnossen-2023-twisted-ambi}.
    We refer the reader to~\cite[Proposition~3.7]{Cnossen-2023-twisted-ambi} for a modern discussion of dualizing objects and the universal property of the norm maps. These constructions were developed in various forms in~\cite{Klein-2001-dualizing-object, Bauer-2004-pcompact, May-Sigurdsson-2006-parametrized, Rognes-2008-Galois, Hopkins-Lurie-2013-ambi, Costenoble-Waner-2016-equivariant, Nikolaus-Scholze-2018-TC, Balmer-DellAmbrogio-2020-Mackey}.
    For the relation with parametrized assembly maps, see~\cite{Quigley-Shah-2025-parametrized-Tate}.
    
    \begin{definition}
        We say that $f$ is \emph{twisted $\cC$-ambidextrous} if the norm is a natural isomorphism.
        If $f$ is truncated, we say it is \emph{$\cC$-ambidextrous} if $f$ and all its iterated diagonals are twisted $\cC$-ambidextrous. In this case, $\DD_f \simeq A^*\ounit_{\cC}$ (\cite[Proposition~3.17]{Cnossen-2023-twisted-ambi}).
    \end{definition}

    \begin{example}[{\cite[Example~4.3.11]{Hopkins-Lurie-2013-ambi},~\cite[Corollary~5.3]{BCSY-Fourier}}]
    \label{exm:ModB-ambi}
        Any truncated anima is $\ModB$-ambidextrous.
    \end{example}
    
    A useful criterion for faithfulness of Galois extensions is given by the norm map. If $R$ is an $A$-Galois extension, it is in particular an $A$-local system. Thus there is a norm map
    \begin{equation*}
        \Nm \;\colon\; A_!(R \otimes \DD_A) \;\too\; A_* R.
    \end{equation*}
    
    Rognes showed that for a connected \COB-dualizable anima $A$ and any $A$-Galois extension $R$, its dual is $R\dual \simeq R \otimes \DD_A$ and that $R$ is faithful if and only if the norm is an isomorphism~\cite[Propositions~6.3.3 and~6.4.7]{Rognes-2008-Galois}. In particular, $R$ is faithful if and only if $R\dual$ satisfies~\labelcref{cond:coGal-1}, which is the obstruction for $R\dual$ to be coGalois.
    We repeat some proofs here for the reader's convenience, and extend the results to non-connected animae.

    \begin{lemma}
    \label{lem:X^G=X@1^G}
        Let $X \in \cC^A$. Then
        \begin{equation*}
            \Delta_* X \;\simeq\; \hom(\ounit_{\cC}[\Path A], \pr_1^* X).
        \end{equation*}
        In particular, if $A$ is \COB-dualizable, then
        \begin{equation*}
            \Delta_* X \;\simeq\; \pr_1^* X \otimes \ounit_{\cC}^{\Path A}.
        \end{equation*}
    \end{lemma}
    \begin{proof}
        Using~\cref{cor:internal-adjoints-maps-of-animae},
        \begin{equation*}
            \begin{split}
                \Delta_* X 
                & \;\simeq\; \Delta_* \hom(A^* \ounit_{\cC}, X) \\
                & \;\simeq\; \Delta_*\hom(A^* \ounit_{\cC}, \Delta^* \pr_1^* X) \\
                & \;\simeq\; \hom(\Delta_! A^* \ounit_{\cC}, \pr_1^* X) \\
                & \;\simeq\; \hom(\ounit_{\cC}[\Path A], \pr_1^* X).
            \end{split}
        \end{equation*}
    \end{proof}


    \begin{corollary}
    \label{cor:M@R=M^G}
        Let $A$ be a \COB-dualizable anima and $R$ a (not necessarily faithful) $A$-Galois extension in $\cC$.
        Then for any $M \in \Mod_R(\cC^A)$,
        \begin{equation*}
            \pr_1^* M \otimes \pr_2^* R \;\simeq\; \Delta_* M \qand M \simeq (\pr_1)_* (\pr_1^* M \otimes \pr_2^* R).
        \end{equation*}
    \end{corollary}
    \begin{proof}
        By~\labelcref{cond:Galois-2} and~\cref{lem:X^G=X@1^G}:
        \begin{equation*}
            \begin{split}
                \pr_1^* M \otimes \pr_2^* R 
                & \;\simeq\; \pr_1^* M \otimes_{\pr_1^* R} \pr_1^* R \otimes \pr_2^* R \\
                & \;\simeq\; \pr_1^* M \otimes_{\pr_1^* R} \pr_1^* R \otimes \ounit_{\cC}^{\Path A} \\
                & \;\simeq\; \pr_1^* M \otimes \ounit_{\cC}^{\Path A} \\
                & \;\simeq\; \Delta_* M.
            \end{split}
        \end{equation*}
        The second part follows as
        \begin{equation*}
            M \;\simeq\; (\pr_1)_* \Delta_* M.
        \end{equation*}
    \end{proof}

    \begin{lemma}
    \label{lem:free-induced-splits}
        Let $X\in \cC^A$. Then the compositions
        \begin{equation*}
            \begin{split}
                & \Delta_! X \;\simeq\; \Delta_!(\Delta^* \pr_1^* X \otimes A^* \ounit_{\cC}) \;\xtoo{\Lproj}\; \pr_1^* X \otimes \Delta_! A^* \ounit_{\cC} \;\simeq\; \pr_1^* X \otimes \ounit_{\cC}[\Path A] \\
                & \Delta_! X \;\simeq\; \Delta_!(\Delta^* \pr_2^* X \otimes A^* \ounit_{\cC}) \;\xtoo{\Lproj}\; \pr_2^* X \otimes \Delta_! A^* \ounit_{\cC} \;\simeq\; \pr_2^* X \otimes \ounit_{\cC}[\Path A]
            \end{split}
        \end{equation*}
        are isomorphisms.
    \end{lemma}
    \begin{proof}
        $\Delta^*$ satisfies the left projection formula by~\cref{exm:f!-projection-formula}.
    \end{proof}

    \begin{corollary}
    \label{cor:hom-R-R}
        Let $A\in \An$ be \COB-dualizable and $R$ be a (not necessarily faithful) $A$-Galois extension.
        Let $X \in \cC^A$. Then $\pr_1^* X \otimes \Delta_! R \simeq \hom(\pr_1^* R, \pr_1^* X \otimes \pr_2^* R)$ in $\cC^{A\times A}$.
    \end{corollary}
    \begin{proof}
        By~\cref{lem:free-induced-splits,cor:internal-adjoints-maps-of-animae} and the \CO-dualizability of $A$,
        \begin{equation*}
            \begin{split}
                \pr_1^* X \otimes \Delta_! R 
                & \;\simeq\; \pr_1^* X \otimes \pr_2^* R \otimes \Delta_! A^* \ounit_{\cC} \\
                & \;\simeq\; \hom(\Delta_* A^* \ounit_{\cC}, \pr_1^* X \otimes \pr_2^* R) \\
                & \;\simeq\; \hom_{\pr_2^* R}(\pr_2^* R \otimes \Delta_* A^* \ounit_{\cC}, \pr_1^* X \otimes \pr_2^* R).
            \end{split}
        \end{equation*}
        Since $A^* \ounit_{\cC}$ is self-dual and $\Delta_! A^* \ounit_{\cC} = \ounit_{\cC}[\Path A]$ is dualizable, by~\cref{lem:Y-and-f!Y-dual-projection} the projection map is an isomorphism
        \begin{equation*}
            \pr_2^* R \otimes \Delta_* A^* \ounit_{\cC} \;\isotoo\; \Delta_*(\Delta^* \pr_2^* R \otimes A^* \ounit_{\cC}) \;\simeq\; \Delta_* R.
        \end{equation*}
        Therefore, by the second Galois condition~\labelcref{cond:Galois-2},
        \begin{equation*}
            \begin{split}
                \pr_1^* X \otimes \Delta_! R 
                & \;\simeq\; \hom_{\pr_2^* R}(\Delta_* R, \pr_1^* X \otimes \pr_2^* R) \\
                & \;\simeq\; \hom_{\pr_2^* R}(\pr_1^* R \otimes \pr_2^* R, \pr_1^* X \otimes \pr_2^* R) \\
                & \;\simeq\; \hom(\pr_1^* R, \pr_1^*X \otimes \pr_2^* R).
            \end{split}
        \end{equation*}
    \end{proof}
    
    \begin{construction}
    \label{cons:comparison-map}
        Let $X\in \cC^A$. There are canonical comparison maps in $\cC^A$
        \begin{equation*}
            \begin{split}
                & \comp \;\colon\; X \otimes \DD_A \;\simeq\; X \otimes (\pr_1)_* \Delta_! A^* \ounit_{\cC} 
                \;\xtoo{\Rproj}\; (\pr_1)_* (\pr_1^* X \otimes \Delta_! A^* \ounit_{\cC})
                \;\isogetss\; (\pr_1)_* \Delta_! X, \\
                & \comp \;\colon\; X \otimes \DD_A \;\simeq\; X \otimes (\pr_2)_* \Delta_! A^* \ounit_{\cC} 
                \;\xtoo{\Rproj}\; (\pr_2)_* (\pr_2^* X \otimes \Delta_! A^* \ounit_{\cC})
                \;\isogetss\; (\pr_2)_* \Delta_! X,
            \end{split}
        \end{equation*}
        where the last isomorphisms are those of~\cref{lem:free-induced-splits}.
    \end{construction}

    \begin{corollary}
    \label{cor:X[G]_hG=X@D}
        If $X\in \cC^A$ is dualizable, then the comparison maps are isomorphisms.
    \end{corollary}
    \begin{proof}
        This follows from~\cref{lem:projection-for-dbl}.
    \end{proof}

    We restrict our attention to dualizable Galois extensions. In~\cref{prop:Gal-dualizable} we show that every Galois extension is dualizable.

    \begin{proposition}[{\cite[Proposition~6.4.7]{Rognes-2008-Galois}}]
    \label{prop:dual-of-Galois}
        Let $A$ be \COB-dualizable and let $R$ be a dualizable (not necessarily faithful) $A$-Galois extension. Then its weak dual is
        \begin{equation*}
            \hom(R, A^* \ounit_{\cC}) \;\simeq\; R \otimes \DD_A.
        \end{equation*}
    \end{proposition}
    \begin{proof}
        As $R$ is dualizable, the comparison map
        \begin{equation*}
            \comp \;\colon\; R \otimes \DD_A \;\too\; (\pr_1)_*\Delta_! R
        \end{equation*}
        is an isomorphism. Using~\cref{cor:hom-R-R,cor:internal-adjoints-maps-of-animae}, the Beck-Chevalley condition, and the first Galois condition~\labelcref{cond:Galois-1}, there is an isomorphism
        \begin{equation*}
            \begin{split}
                (\pr_1)_* \Delta_! R
                & \;\simeq\; (\pr_1)_* \hom(\pr_1^* R, \pr_2^* R) \\
                & \;\simeq\; \hom(R, (\pr_1)_* \pr_2^* R) \\
                & \;\simeq\; \hom(R, A^* A_* R) \\
                & \;\simeq\; \hom(R, A^* \ounit_{\cC}).
            \end{split}
        \end{equation*}
    \end{proof}

    \begin{corollary}
    \label{cor:dual-of-trivial-Galois}
        Let $A$ be \COB-dualizable. Then the weak dual of the universal coGalois extension is
        \begin{equation*}
            \ounit_{\cC}^{\Path A} \;\simeq\; \ounit_{\cC}[\Path A] \otimes \DD_{\pr_1} \qin \cC^{A \times A}.
        \end{equation*}
    \end{corollary}
    \begin{proof}
        We first show that the local system of algebras $\ounit_{\cC}^{\Path A} \in \CAlg(\cC^A)^A$ is Galois. 
        We use the coGalois conditions for its dual $\ounit_{\cC}[\Path A]$.
        Since $A$ is \COB-dualizable, by~\cref{cor:internal-adjoints-maps-of-animae},
        \begin{equation*}
            \ounit_{\cC}^{\Path A} \boxtimes \ounit_{\cC}^{\Path A} \;\isotoo\; \Delta_* \ounit_{\cC}^{\Path A} \qin \CAlg(\cC^A)^{A \times A}.
        \end{equation*}
        Thus it satisfies the second Galois condition.
        
        For the first Galois condition, since $(-)\dual \simeq \hom(-, A^* \ounit_{\cC})$, it sends colimits in $\cC^A$ to limits in $\cC^A$. In particular,
        \begin{equation*}
            A_* \ounit_{\cC}^{\Path A} \;=\; A_* (\ounit_{\cC}[\Path A]\dual) \;\simeq\; (A_! \ounit_{\cC}[\Path A])\dual \;\isogetss\; (A^*\ounit_{\cC})\dual \;=\; A^*\ounit_{\cC},
        \end{equation*}
        where the map is the dual of the counit of $\ounit_{\cC}[\Path A]$, i.e.\ the unit of $\ounit_{\cC}^{\Path A}$.

        The local system $\ounit_{\cC}^{\Path A}$ is dualizable, as the dual of the dualizable local system $\ounit_{\cC}[\Path A]$.
        The claim now follows from~\cref{prop:dual-of-Galois}, using that the dualizing object of $A$ in $\cC$ identifies with the dualizing object of $\pr_1$ in $\cC^A$.
    \end{proof}

    \begin{lemma}
    \label{lem:induced-tate-is-zero}
        Let $A$ be \COB-dualizable.
        \begin{enumerate}[label=(\arabic*)]
            \item \label{item:induced-tate-pr1}
            For any $X, Y \in \cC^A$, the norm map
            \begin{equation*}
                \Nm \colon (\pr_1)_! (\Delta_! X \otimes \pr_1^* Y \otimes \DD_{\pr_1}) \;\too\; (\pr_1)_* (\Delta_! X \otimes \pr_1^* Y) \qin \cC^{A}
            \end{equation*}
            is an isomorphism.
            \item \label{item:induced-tate-pr2}
            For any dualizable $X \in \cC^A$ and any $Y \in \cC^A$, the norm map
            \begin{equation*}
                \Nm \colon (\pr_1)_! (\Delta_* X \otimes \pr_2^* Y \otimes \DD_{\pr_1}) \;\too\; (\pr_1)_* (\Delta_* X \otimes \pr_2^* Y) \qin \cC^{A}
            \end{equation*}
            is an isomorphism.
        \end{enumerate}
    \end{lemma}
    \begin{proof}
        For~\labelcref{item:induced-tate-pr1}, the projection formula of~\cref{exm:f!-projection-formula} identifies the source with
        \begin{equation*}
            (\pr_1)_! (\Delta_! X \otimes \pr_1^* Y \otimes \DD_{\pr_1}) 
            \;\simeq\; (\pr_1)_! \Delta_!(X \otimes \Delta^* (\pr_1^* Y \otimes \DD_{\pr_1})) 
            \;\simeq\; X \otimes Y \otimes \Delta^* \DD_{\pr_1}.
        \end{equation*}
        By~\cref{lem:free-induced-splits,cor:dual-of-trivial-Galois}, the target identifies as
        \begin{equation*}
            (\pr_1)_* (\Delta_! X \otimes \pr_1^* Y)
            \;\simeq\; (\pr_1)_* (\pr_1^* X \otimes \ounit_{\cC}[\Path A] \otimes \pr_1^* Y) 
            \;\simeq\; (\pr_1)_* (\pr_1^* (X \otimes Y)  \otimes \ounit_{\cC}^{\Path A} \otimes \DD_{\pr_1}).
        \end{equation*}
        Now, using the dualizability of $\ounit_{\cC}[\Path A]$,
        \begin{equation*}
            \begin{split}
                (\pr_1)_* (\Delta_! X \otimes \pr_1^* Y) 
                & \;\simeq\; (\pr_1)_* \hom(\ounit_{\cC}[\Path A], \pr_1^* (X \otimes Y) \otimes \DD_{\pr_1}) \\
                & \;\simeq\; (\pr_1)_* \hom(\Delta_! A^* \ounit_{\cC}, \pr_1^* (X \otimes Y) \otimes \DD_{\pr_1}).
            \end{split}
        \end{equation*}        
        By~\cref{cor:internal-adjoints-maps-of-animae},
        \begin{equation*}
            (\pr_1)_* (\Delta_! X \otimes \pr_1^* Y) 
            \;\simeq\; (\pr_1)_* \Delta_* \hom(A^* \ounit_{\cC}, \Delta^* \pr_1^* (X \otimes Y) \otimes \Delta^* \DD_{\pr_1}) 
            \;\simeq\; X \otimes Y \otimes \Delta^* \DD_{\pr_1}.
        \end{equation*}
        Tracing the norm map through these identifications shows that it is an isomorphism.

        We now prove~\labelcref{item:induced-tate-pr2}. Since $X$ is dualizable, \cref{lem:free-induced-splits} gives
        \begin{equation*}
            \Delta_!(X\dual) \;\simeq\; \pr_1^*(X\dual)\otimes\ounit_{\cC}[\Path A],
        \end{equation*}
        which is dualizable by the \COB-dualizability of $A$. Hence, by~\cref{lem:Y-and-f!Y-dual-projection} with the map $\Delta$,
        \begin{equation*}
            \Delta_*X \otimes \pr_2^*Y \;\simeq\; \Delta_*(X \otimes Y).
        \end{equation*}
        Moreover, \cref{lem:X^G=X@1^G}, \cref{cor:dual-of-trivial-Galois} and~\cref{lem:free-induced-splits} give
        \begin{equation*}
            \Delta_*(X\otimes Y)\otimes\DD_{\pr_1} \;\simeq\; \pr_1^*(X\otimes Y)\otimes\ounit_{\cC}^{\Path A}\otimes\DD_{\pr_1} \;\simeq\; \pr_1^*(X\otimes Y)\otimes\ounit_{\cC}[\Path A] \;\simeq\; \Delta_!(X\otimes Y).
        \end{equation*}
        It follows that the source and target of the norm map identify as
        \begin{equation*}
            \begin{split}
                & (\pr_1)_!(\Delta_*X\otimes\pr_2^*Y\otimes\DD_{\pr_1}) \;\simeq\; (\pr_1)_!\Delta_!(X\otimes Y) \;\simeq\; X\otimes Y, \\
                & (\pr_1)_*(\Delta_*X\otimes\pr_2^*Y) \;\simeq\; (\pr_1)_*\Delta_*(X\otimes Y) \;\simeq\; X\otimes Y.
            \end{split}
        \end{equation*}
        Tracing the norm map through these identifications shows that it is the identity.
    \end{proof}

    \begin{corollary}
    \label{cor:Galois-norm-is-iso}
        Let $A$ be \COB-dualizable and let $R\in\Gal_A(\cC)$ be dualizable. Then for every $Y \in \cC^A$, the norm map
        \begin{equation*}
            \Nm \;\colon\; A_!(R\otimes Y\otimes\DD_A) \;\too\; A_*(R\otimes Y)
        \end{equation*}
        is an isomorphism.
    \end{corollary}
    \begin{proof}
        Since $R$ is faithful, it is enough to show the norm is an isomorphism after applying the conservative functor $R\otimes A^*(-)\colon\cC\to\cC^A$.
        By Beck--Chevalley, the left projection formula, the second Galois condition~\labelcref{cond:Galois-2}, and the identification $\pr_2^*\DD_A\simeq\DD_{\pr_1}$, its source becomes
        \begin{equation*}
            \begin{split}
                R\otimes A^*A_!(R\otimes Y \otimes\DD_A) 
                & \;\simeq\; R \otimes (\pr_1)_! \pr_2^*(R \otimes Y \otimes \DD_A) \\
                & \;\simeq\; (\pr_1)_!(\pr_1^*R \otimes \pr_2^*R \otimes\pr_2^*Y \otimes \pr_2^*\DD_A) \\
                & \;\simeq\; (\pr_1)_!(\Delta_*R \otimes \pr_2^*Y \otimes \DD_{\pr_1}).
            \end{split}
        \end{equation*}
        Similarly, Beck--Chevalley, the right projection formula for the dualizable object $R$, and~\labelcref{cond:Galois-2} identify its target as
        \begin{equation*}
            \begin{split}
                R\otimes A^*A_*(R \otimes Y) 
                & \;\simeq\; R \otimes (\pr_1)_*\pr_2^*(R \otimes Y) \\
                & \;\simeq\; (\pr_1)_*(\pr_1^*R \otimes \pr_2^*R \otimes \pr_2^*Y) \\
                & \;\simeq\; (\pr_1)_*(\Delta_*R \otimes \pr_2^*Y).
            \end{split}
        \end{equation*}
        Under these identifications, the resulting map is the norm map of~\cref{lem:induced-tate-is-zero}\labelcref{item:induced-tate-pr2}, with $X=R$, and is therefore an isomorphism.
    \end{proof}

    \begin{proposition}[{\cite[Proposition~6.2.1]{Rognes-2008-Galois}}]
    \label{prop:Gal-dualizable}
        Let $A \in \An$ be \COB-dualizable and $R$ a (not necessarily faithful) $A$-Galois extension. Then $R\in \cC^A$ is dualizable.
    \end{proposition}
    \begin{proof}
        It is enough to show that
        \begin{equation*}
            M \otimes \hom(R, A^*\ounit_{\cC}) \;\isotoo\; \hom(R, M)
        \end{equation*}
        for any $M \in \Mod_R(\cC^A)$.
    
        This follows from the commutative diagram:
        \begin{equation*}
            \begin{tikzcd}
            	{M \otimes \hom(R, A^*\ounit_{\cC})} & {\hom(R,M)} \\
            	{M \otimes \hom(R, (\pr_1)_* \pr_2^* R)} & {\hom(R, (\pr_1)_*(\pr_1^* M \otimes \pr_2^* R))} \\
            	{M \otimes (\pr_1)_*\hom(\pr_1^* R, \pr_2^* R)} & {(\pr_1)_* \hom(\pr_1^* R, \pr_1^* M \otimes \pr_2^* R)} \\
            	{M \otimes (\pr_1)_* \Delta_! R} & {(\pr_1)_* (\pr_1^* M \otimes \Delta_! R)} \\
            	{M \otimes (\pr_1)_! (\Delta_!R \otimes \DD_{\pr_1})} & {(\pr_1)_!(\pr_1^* M \otimes \Delta_! R \otimes \DD_{\pr_1}).}
            	\arrow[from=1-1, to=1-2]
            	\arrow["{\labelcref{cor:M@R=M^G}}", from=1-2, to=2-2]
            	\arrow["\begin{array}{c} \begin{array}{c} \begin{array}{c} \begin{array}{c} \begin{array}{c} \substack{\sim \\ \\} \end{array} \end{array} \end{array} \end{array} \end{array}"{marking, allow upside down}, tail reversed, no head, from=2-1, to=1-1]
                \arrow["{\begin{array}{l}\labelcref{cond:Galois-1} \\ \text{Beck--Chevalley}\end{array}}", , from=1-1, to=2-1]
            	\arrow["{\labelcref{cor:internal-adjoints-maps-of-animae}}", from=2-1, to=3-1]
            	\arrow["\begin{array}{c} \begin{array}{c} \begin{array}{c} \begin{array}{c} \begin{array}{c} \substack{\sim \\ \\} \end{array} \end{array} \end{array} \end{array} \end{array}"{marking, allow upside down}, tail reversed, no head, from=2-2, to=1-2]
            	\arrow["{\labelcref{cor:internal-adjoints-maps-of-animae}}", from=2-2, to=3-2]
            	\arrow["\begin{array}{c} \begin{array}{c} \begin{array}{c} \begin{array}{c} \begin{array}{c} \substack{\sim \\ \\} \end{array} \end{array} \end{array} \end{array} \end{array}"{marking, allow upside down}, tail reversed, no head, from=3-1, to=2-1]
            	\arrow[from=3-1, to=3-2]
            	\arrow["{{{{\labelcref{cor:hom-R-R}}}}}", from=3-1, to=4-1]
            	\arrow["\begin{array}{c} \begin{array}{c} \begin{array}{c} \begin{array}{c} \begin{array}{c} \substack{\sim \\ \\} \end{array} \end{array} \end{array} \end{array} \end{array}"{marking, allow upside down}, tail reversed, no head, from=3-2, to=2-2]
            	\arrow["{{{{\labelcref{cor:hom-R-R}}}}}", from=3-2, to=4-2]
            	\arrow["\begin{array}{c} \begin{array}{c} \begin{array}{c} \begin{array}{c} \begin{array}{c} \substack{\sim \\ \\} \end{array} \end{array} \end{array} \end{array} \end{array}"{marking, allow upside down}, tail reversed, no head, from=4-1, to=3-1]
            	\arrow[from=4-1, to=4-2]
            	\arrow["\begin{array}{c} \begin{array}{c} \begin{array}{c} \begin{array}{c} \begin{array}{c} \substack{\sim \\ \\} \end{array} \end{array} \end{array} \end{array} \end{array}"{marking, allow upside down}, tail reversed, no head, from=4-2, to=3-2]
            	\arrow["\begin{array}{c} \begin{array}{c} \begin{array}{c} \begin{array}{c} \begin{array}{c} \substack{\sim \\ \\} \end{array} \end{array} \end{array} \end{array} \end{array}"{marking, allow upside down}, from=5-1, to=4-1]
            	\arrow["{{{{\labelcref{lem:induced-tate-is-zero}}}}}"', from=5-1, to=4-1]
            	\arrow["\sim", from=5-1, to=5-2]
            	\arrow["{{\labelcref{exm:f!-projection-formula}}}"', from=5-1, to=5-2]
            	\arrow["{{{{\labelcref{lem:induced-tate-is-zero}}}}}"', from=5-2, to=4-2]
            	\arrow["\begin{array}{c} \begin{array}{c} \begin{array}{c} \begin{array}{c} \begin{array}{c} \substack{\sim \\ \\} \end{array} \end{array} \end{array} \end{array} \end{array}"{marking, allow upside down}, from=5-2, to=4-2]
            \end{tikzcd}
        \end{equation*}
        All vertical maps and the bottom horizontal map are isomorphisms, so the top map is an isomorphism.
    \end{proof}

    \subsubsection*{Faithfulness and the norm map}

    \begin{proposition}[{\cite[Proposition~6.3.3]{Rognes-2008-Galois}}]
    \label{prop:faithful<=>Nm-is-iso}
        Let $A$ be a \COB-dualizable anima and $R$ a (not necessarily faithful) $A$-Galois extension in $\cC$. Then $R$ is faithful if and only if the norm map
        \begin{equation*}
            \Nm \;\colon\; A_!(R \otimes \DD_A) \;\too\; A_* R \;\isogetss\; \ounit_{\cC}
        \end{equation*}
        is an isomorphism.
    \end{proposition}
    \begin{proof}
        Assume that the norm map is an isomorphism. Let $f \colon X \to Y$ in $\cC$ be such that
        \begin{equation*}
            A^* f \otimes \id \;\colon\; A^* X \otimes R \;\too\; A^* Y \otimes R \qin \cC^A
        \end{equation*}
        is an isomorphism. Consider the commutative diagram
        \begin{equation*}
            \begin{tikzcd}
                {A_!(A^* X \otimes R\otimes \DD_A)} && {A_!(A^* Y \otimes R \otimes \DD_A)} \\
                {X \otimes A_!(R \otimes \DD_A)} && {Y \otimes A_!(R \otimes \DD_A)} \\
                {X \otimes A_*R} && {Y \otimes A_*R} \\
                X && Y.
                \arrow["{A_!(A^* f \otimes \id)}", from=1-1, to=1-3]
                \arrow["\Lproj", from=1-1, to=2-1]
                \arrow["\Lproj", from=1-3, to=2-3]
            	\arrow["{f \otimes \id}", from=2-1, to=2-3]
                \arrow["\Nm", from=2-1, to=3-1]
                \arrow["\Nm", from=2-3, to=3-3]
                \arrow["{f \otimes \id}", from=3-1, to=3-3]
                \arrow["{\rotatebox{90}{$\sim$}}"', from=4-1, to=3-1]
                \arrow["f", from=4-1, to=4-3]
                \arrow["{\rotatebox{90}{$\sim$}}"', from=4-3, to=3-3]
            \end{tikzcd}
        \end{equation*}
        By assumption, the top horizontal map and the norm maps are invertible. The projection maps are invertible by~\cref{exm:f!-projection-formula}. Therefore $f$ is invertible, as needed.

        Conversely, assume that $R$ is faithful. Then $R$ is dualizable by~\cref{prop:Gal-dualizable}, and the norm map is an isomorphism by~\cref{cor:Galois-norm-is-iso}, applied to $Y = A^* \ounit_{\cC}$.
    \end{proof}


    We are finally ready to compare Galois and coGalois extensions:
    
    \begin{proposition}
    \label{prop:Gal-coGal-iso}
        Assume $A$ is \COB-dualizable. Then the dual functor induces an isomorphism of animae
        \begin{equation*}
            (-)\dual \;\colon\; \Gal_A(\cC) \;\isotoo\; \coGal_A(\cC).
        \end{equation*}
    \end{proposition}
    \begin{proof}
        Note that any coGalois extension is dualizable by~\cref{cor:coGal-dualizable} and any faithful Galois extension (and indeed any Galois extension) is dualizable by~\cref{prop:Gal-dualizable}.
        
        The dual defines an isomorphism between dualizable $A$-local systems of algebras and dualizable $A$-local systems of coalgebras:
        \begin{equation*}
            (-)\dual \;\colon\; (\CAlg(\cC)^A)^{\mrm{dbl},\simeq} \;\isotoo\; (\coCAlg(\cC)^A)^{\mrm{dbl},\simeq}.
        \end{equation*}
        Let $R\in \CAlg(\cC)^A$ be dualizable, with $Q \coloneq R\dual$. Since $A$ is \COB-dualizable, by~\cref{lem:free-induced-splits,cor:internal-adjoints-maps-of-animae}, $R$ satisfies the second Galois condition
        \begin{equation*}
            R \boxtimes R \;\isotoo\; \Delta_* R
        \end{equation*}
        if and only if $Q$ satisfies the second coGalois condition
        \begin{equation*}
            \Delta_! Q \;\isotoo\; Q \boxtimes Q.
        \end{equation*}
        It is therefore left to show that the dual of a faithful Galois extension satisfies the first coGalois condition, and that the dual of a coGalois extension is faithful and satisfies the first Galois condition.
        
        Let $R$ be a (dualizable) $A$-Galois extension with dual $Q = R\dual$. Tracing the maps in the proof of~\cref{prop:dual-of-Galois}, under the isomorphism $Q \simeq R \otimes \DD_A$, the counit $A_! Q \to \ounit_{\cC}$ is identified with the norm map. Thus, by~\cref{prop:faithful<=>Nm-is-iso}, $R$ is faithful if and only if $Q$ satisfies the first coGalois condition. In particular, the dual of a faithful Galois extension is a coGalois extension.

        On the other hand, if $Q$ is a coGalois extension with dual $R = Q\dual$, then since $(-)\dual \simeq \hom(-, A^*\ounit_{\cC})$, it sends colimits to limits. In particular, by~\cref{cor:internal-adjoints-maps-of-animae},
        \begin{equation*}
            A_* R \;=\; A_* (Q\dual) \;\simeq\; (A_! Q)\dual \;\isogetss\; \ounit_{\cC}\dual \;=\; \ounit_{\cC},
        \end{equation*}
        where the map is the dual of the counit of $Q$, i.e.\ the unit of $R$. Thus $R$ is a dualizable $A$-Galois extension. As $Q$ is coGalois, $R$ is faithful by~\cref{prop:faithful<=>Nm-is-iso}.
    \end{proof}

    \begin{theorem}
    \label{thm:Galois-corepresented}
        Assume $A$ is \COB-dualizable. 
        Then faithful $A$-Galois extensions are corepresented: there is an isomorphism
        \begin{equation*}
            \begin{split}
                \MapBo(\basecat^A, \cC) & \;\isotoo\; \Gal_A(\cC) \\
                F & \;\mapsto\; F^A(\ounit_{\basecat}[\Path A])\dual \;\simeq\; F^A(\ounit_{\basecat}^{\Path A}).
            \end{split}
        \end{equation*}
        Its inverse sends $R \in \Gal_A(\cC)$ to the functor
        \begin{equation*}
            F_R \;\colon\; \basecat^A \;\xtoo{\ounit_{\cC}[-]}\; \cC^A \;\xtoo{R\otimes -}\; \cC^A \;\xtoo{A_*}\; \cC.
        \end{equation*}
    \end{theorem}
    \begin{proof}
        The isomorphism follows from~\cref{prop:coGal-FunL-iso} and~\cref{prop:Gal-coGal-iso}. Its inverse sends $R \in \Gal_A(\cC)$ to the functor
        \begin{equation*}
            F_R(X) \;=\; A_!(R\dual \otimes \ounit_{\cC}[X]) \;\simeq\; A_!(R \otimes \ounit_{\cC}[X] \otimes \DD_A),
        \end{equation*}
        where the second isomorphism is by~\cref{prop:dual-of-Galois}.
        By~\cref{cor:Galois-norm-is-iso,prop:Gal-dualizable}, the norm gives an isomorphism
        \begin{equation*}
            F_R(X) \;\simeq\; A_!(R \otimes \ounit_{\cC}[X] \otimes \DD_A) \;\xisotoo{\Nm}\; A_*(R \otimes \ounit_{\cC}[X]),
        \end{equation*}
        which proves the formula for the inverse.
    \end{proof}

    \begin{example}
    \label{exm:tensor-of-two-galois-extensions}
        Let $A,B$ be two \COB-dualizable animae and let $R \in \Gal_A(\cC)$, $S \in \Gal_B(\cC)$. Then $A \times B$ is again \COB-dualizable and $R \boxtimes S \in \Gal_{A \times B}(\cC)$. Its corresponding $\basecat$-linear symmetric monoidal functor is given by
        \begin{equation*}
            F_{R \boxtimes S} \;\colon\; \basecat^{A \times B} \;\simeq\; \basecat^A \otimes_{\basecat} \basecat^B \;\xtoo{F_R \otimes F_S} \cC \otimes_{\basecat} \cC \;\xtoo{\otimes}\; \cC,
        \end{equation*}
        where the last map is the tensor product in $\cC$.
    \end{example}

\section{Galois theory and affineness}
\label{sec:affine-galois-theory}
We develop some of the cornerstones of Galois theory in our generality. 
We first discuss trivial Galois extensions in~\cref{subsec:trivial-extensions}. 
In \cref{subsec:affineness}, we recall the notion of affineness, as introduced in~\cite{BCSY-Fourier} (following~\cite{Gaitsgory-2015-1-affineness,Mathew-Meier-2015-affineness}), and study its interaction with dualizability and categorification. 
We then prove Galois descent in~\cref{subsec:Galois-descent}. 
Finally, in~\cref{subsec:relative-Galois}, we study relative Galois theory, the analogue of the interaction of quotients and subgroups with Galois theory.
The role of quotients will be played by general maps of appropriate animae, while we replace the notion of subgroups by that of affine maps.

We begin by showing that Galois extensions are functorial:
\begin{proposition}
\label{cor:Galois-preserved-under-symmetric-monoidal}
    Let $F\colon \cC \to \cD$ be a map in $\CAlgB$. Then, for every \COB-dualizable anima $A$, the functor $F^A$ sends faithful $A$-Galois extensions in $\cC$ to faithful $A$-Galois extensions in $\cD$.
\end{proposition}

\begin{proof}
    Let $A$ be \COB-dualizable and let $R\in \Gal_A(\cC)$. By~\cref{thm:Galois-corepresented}, $R$ corresponds to a $\basecat$-linear symmetric monoidal functor $F_R\colon \basecat^A \to \cC$. The composition $\basecat^A \xto{F_R} \cC \xto{F} \cD$ is symmetric monoidal and therefore corresponds to the faithful $A$-Galois extension
    \begin{equation*}
        F^A(F_R^A(\ounit_{\basecat}^{\Path A})) \;\simeq\; F^A(R).
    \end{equation*}
\end{proof}

\subsection{Trivial Galois extensions}
\label{subsec:trivial-extensions}

    Every point $a \in A$ gives a symmetric monoidal evaluation functor $a^* \colon \cC^A \to \cC$. If $A$ is \CO-dualizable, this determines a distinguished $A$-Galois extension by corepresentability.

    \begin{definition}
    \label{def:trivial-Galois}
        Let $A$ be a \CO-dualizable anima and $a \in A$. We define the trivial $A$-Galois extension $\triv_a =\triv_{A,a}$ as corresponding to the symmetric monoidal functor
        \begin{equation*}
            \ev_a = a^* \;\colon\; \cC^A \;\too\; \cC.
        \end{equation*}
    \end{definition}

    \begin{remark}
        Equivalently, $\triv_a$ is the $A$-local system of algebras whose value at $b \in A$ is
        \begin{equation*}
            (\triv_a)_b \;\simeq\; \ounit_{\cC}^{\Path_{ab} A} \;\simeq\; \begin{cases}
                \ounit_{\cC}^{\Omega_a A}, & a \sim b, \\
                0, & a \not\sim b.
            \end{cases}
        \end{equation*}
    \end{remark}

    \begin{example}
        For $A = \BG$ with a base point $e \in \BG$, the trivial Galois extension $\triv_{\BG,e}$ is the usual split $G$-Galois extension $\ounit_{\cC}^G$; see~\cite[\S5.1]{Rognes-2008-Galois}.
    \end{example}
    
    We now package the trivial extensions into an $A$-local system of symmetric monoidal functors.

    \begin{definition}
        Let $\cC \in \CAlg(\Pr)$ and $A \in \An$. The lax monoidality of the functor 
        \begin{equation*}
            \Delta_* \;\colon\; (\ModC)^A \;\too\; (\ModC)^{A \times A}
        \end{equation*}
        endows $\Delta_* A^* \cC$ with an algebra structure. Using that 
        \begin{equation*}
            (\pr_2)_* \;\colon\; (\ModC)^{A \times A} \;\too\; (\ModC)^A
        \end{equation*}
        is also lax symmetric monoidal, and applying it to the unit map of $\Delta_* A^* \cC$, we get a symmetric monoidal functor
        \begin{equation*}
            \ev_A \;\colon\; A^* \cC^A \;\simeq\; (\pr_2)_* (A \times A)^* \cC \;\too\; (\pr_2)_* \Delta_* A^* \cC \;\simeq\; A^* \cC \qin \CAlgC^A.
        \end{equation*}
    \end{definition}

    \begin{lemma}
    \label{lem:trivA-at-a-point}
        Let $\cC \in \CAlg(\Pr)$, $A \in \An$ and $a \in A$.
        Then the functor
        \begin{equation*}
            \cC^A \;\simeq\; a^* A^* \cC^A \;\xtoo{a^* \ev_A}\; a^* A^* \cC \;\simeq\; \cC
        \end{equation*}
        identifies with $a^*$.
    \end{lemma}
    \begin{proof}
        Consider the pullback square
        \begin{equation*}
            \begin{tikzcd}
                A & \pt \\
                {A \times A} & A.
                \arrow[from=1-1, to=1-2]
                \arrow["{(\id, a)}"', from=1-1, to=2-1]
                \arrow["\lrcorner"{anchor=center, pos=0.125}, draw=none, from=1-1, to=2-2]
                \arrow["a", from=1-2, to=2-2]
                \arrow["{\pr_2}", from=2-1, to=2-2]
            \end{tikzcd}
        \end{equation*}
        By the Beck--Chevalley condition, the following diagram commutes:
        \begin{equation*}
            \begin{tikzcd}
                {(\ModC)^{A\times A}} && {(\ModC)^A} \\
                \\
                {(\ModC)^A} && (\ModC).
                \arrow["{(\pr_2)_*}", from=1-1, to=1-3]
                \arrow["{(\id, a)^*}"', from=1-1, to=3-1]
                \arrow["{a^*}", from=1-3, to=3-3]
                \arrow["{A_*}", from=3-1, to=3-3]
            \end{tikzcd}
        \end{equation*}
        The unit map $(A \times A)^* \cC \to \Delta_* A^* \cC$ in $(\ModC)^{A \times A}$ is sent under the top-right composition to $a^* \ev_A$. Its image under $(\id, a)^*$ is the unit map
        \begin{equation*}
            A^* \cC \;\too\; (\id, a)^* \Delta_* A^* \cC.
        \end{equation*}
        Using the Beck--Chevalley condition for the pullback diagram
        \begin{equation*}
            \begin{tikzcd}
                \pt & A \\
                A & {A\times A,}
                \arrow["a", from=1-1, to=1-2]
                \arrow["a"', from=1-1, to=2-1]
                \arrow["\lrcorner"{anchor=center, pos=0.125}, draw=none, from=1-1, to=2-2]
                \arrow["{(\id, a)}", from=1-2, to=2-2]
                \arrow["\Delta", from=2-1, to=2-2]
            \end{tikzcd}
        \end{equation*}
        we find that it is the unit map
        \begin{equation*}
            A^* \cC \;\too\; a_* a^* A^* \cC \;\simeq\; a_* \cC.
        \end{equation*}
        Finally, the image of this map under $A_*$ is the map
        \begin{equation*}
            a^* \;\colon\; \cC^A \;\simeq\; A_* A^* \cC \;\too\; A_* a_* \cC \;\simeq\; \cC.
        \end{equation*}
    \end{proof}
    
    \begin{remark}
        If $A$ is \CO-dualizable, the pullback of $\ev_A$ along any $a \in A$ corresponds to the trivial Galois extension $\triv_{A,a}$ of~\cref{def:trivial-Galois}.

        Corepresentability identifies this family with a map
        \begin{equation*}
            \triv \;\colon\; A \;\too\; \Gal_A(\cC), \qquad a \;\mapsto\; \triv_{A,a}.
        \end{equation*}
    \end{remark}

\subsection{Affine maps}
\label{subsec:affineness}
    A $\cC$-affine anima is an anima $A$ such that the category $\cC^A$ identifies symmetric monoidally with $\Mod_{\ounit_{\cC}^A}(\cC)$. This is a strong form of unipotence and plays an important role in the theory. In particular, for such animae the corepresentability theorem can be expressed in terms of maps of commutative algebras. 

    \begin{definition}
        A lax symmetric monoidal functor $f_* \colon \cC \to \cD$ is said to be unital if the unit map $\ounit_{\cD} \to f_*(\ounit_{\cC})$ is an isomorphism.
    \end{definition}
    
    Let $f_* \colon \cC \to \cD$ be lax symmetric monoidal. Then $f_*$ sends algebras in $\cC$ to algebras in $\cD$ and modules to modules. In particular, it induces a unital lax symmetric monoidal functor
    \begin{equation*}
        f_{\scriptscriptstyle\bigstar} \;\colon\; \cC \;\simeq\; \Mod_{\ounit_{\cC}}(\cC) \;\too\; \Mod_{f_*(\ounit_{\cC})} (\cD).
    \end{equation*}
    
    \begin{definition}[{\cite[Definition~2.1]{BCSY-Fourier}}]
        A symmetric monoidal left adjoint functor $f^* \colon \cD \to \cC$ in $\CAlg(\Pr)$, with right adjoint $f_*$, is said to be affine if
        \begin{equation*}
            f_{\scriptscriptstyle\bigstar} \;\colon\; \cC \;\too\; \Mod_{f_*\ounit_{\cC}} (\cD)
        \end{equation*}
        is an isomorphism.
    \end{definition}
    
    \begin{definition}
        Let $\cC \in\CAlg(\Pr)$. A map of animae $f \colon A \to B$ is called $\cC$-affine if $f^* \colon \cC^B \to \cC^A$ is affine.
        We say that $A$ is $\cC$-affine if the terminal map $A \to \pt$ is.
    \end{definition}
    
    We can immediately give an alternative proof of a weaker version of~\cite[Proposition~2.30]{BCSY-Fourier}, assuming \CO-dualizability:
    \begin{corollary}
    \label{cor:corepresntatbility-for-affines}
        Let $\cC \in \CAlg(\Pr)$ and $A$ be a \CO-dualizable anima. Assume $A$ is $\cC$-affine. Then
        \begin{equation*}
            \Gal_A(\cC) \;\simeq\; \Map_{\CAlg(\cC)}(\ounit_{\cC}^A, \ounit_{\cC}).
        \end{equation*}
    \end{corollary}
    \begin{proof}
        By~\cref{thm:Galois-corepresented},
        \begin{equation*}
            \Gal_A(\cC) \;\simeq\; \MapCo(\cC^A, \cC).
        \end{equation*}
        We use the adjunction
        \begin{equation*}
            \Mod_{(-)}(\cC) \;\colon\; \CAlg(\cC) \;\longadj\; \CAlgC \;\cocolon\; \ounit_{(-)}.
        \end{equation*}
        The right adjoint $\ounit_{(-)}$ sends $\cD$ to the $\cC$-valued endomorphism algebra of its unit.
        The result follows using that $\cC^A \simeq \Mod_{\ounit_{\cC}^A}(\cC)$.
    \end{proof}
    
    \begin{lemma}
    \label{lem:affiness-of-ModR}
        Let $\cC \in \CAlg(\Pr)$ and $A \in \An$ be $\cC$-affine. Then for any $R \in \CAlg(\cC)$, $A$ is $\Mod_R(\cC)$-affine.
    \end{lemma}
    \begin{proof}
        Since $A$ is $\cC$-affine,~\cite[Proposition~2.12]{BCSY-Fourier}, applied to $A^*R \in \CAlg(\cC^A)$, gives a symmetric monoidal equivalence
        \begin{equation*}
            \Mod_R(\cC)^A
            \;\simeq\;
            \Mod_{A^*R}(\cC^A)
            \;\xisotoo{A_{\scriptscriptstyle\bigstar}}\;
            \Mod_{A_*A^*R}(\cC)
            \;\simeq\;
            \Mod_{A_*A^*R}(\Mod_R(\cC)).
        \end{equation*}
        This is precisely the affineness of $A$ in $\Mod_R(\cC)$.
    \end{proof}

    \begin{lemma}
    \label{lem:a*-is-affine}
        Let $\cC \in \CAlg(\Pr)$ and $A$ be a \CO-dualizable anima. The map $a \colon \pt \to A$ is $\cC$-affine for all $a \in A$.
    \end{lemma}
    \begin{proof}
        By~\cite[Proposition~2.5]{BCSY-Fourier}, it is enough to show that $a_*$ preserves colimits, is $\cC^A$-linear and is conservative.
        First we show that $a_! \ounit_{\cC}$ is dualizable.
        Using the Beck--Chevalley condition for the pullback diagram
        \begin{equation*}
            \begin{tikzcd}
                \pt & A \\
                A & {A\times A,}
                \arrow["a", from=1-1, to=1-2]
                \arrow["a"', from=1-1, to=2-1]
                \arrow["\lrcorner"{anchor=center, pos=0.125}, draw=none, from=1-1, to=2-2]
                \arrow["{(\id, a)}", from=1-2, to=2-2]
                \arrow["\Delta", from=2-1, to=2-2]
            \end{tikzcd}
        \end{equation*}
        we get that $a_! \ounit_{\cC} \simeq (\id,a)^* \ounit_{\cC}[\Path A] \in \cC^A$, which is dualizable by the \CO-dualizability of $A$. Therefore, by~\cref{cor:internal-adjoints-maps-of-animae},
        \begin{equation*}
            (a_! \ounit_{\cC})\dual \;\simeq\; a_* \ounit_{\cC}.
        \end{equation*}
        
        Again using~\cref{cor:internal-adjoints-maps-of-animae}, for every $Y \in \cC^A$,
        \begin{equation*}
            a_* a^* Y \;\simeq\; a_*\hom(\ounit_{\cC}, a^*Y) \;\simeq\; \hom(a_!\ounit_{\cC}, Y) \;\simeq\; a_*\ounit_{\cC} \otimes Y.
        \end{equation*}
        In particular, for $X \in \cC$
        \begin{equation*}
            a_* X \;\simeq\; a_* a^* A^* X \;\simeq\; a_* \ounit_{\cC} \otimes A^* X,
        \end{equation*}
        and $a_*$ preserves colimits.
        Moreover, for $X \in \cC$ and $Y \in \cC^A$, the projection map is identified with the composition
        \begin{equation*}
            Y \otimes a_*X \;\simeq\; a_*\ounit_{\cC} \otimes Y \otimes A^* X \;\simeq\; a_*a^*(Y \otimes A^*X) \;\simeq\; a_*(a^*Y \otimes X),
        \end{equation*}
        and is therefore an isomorphism. Hence $a_*$ is $\cC^A$-linear. Finally, for $X \in \cC$,
        \begin{equation*}
            a^* a_*X \;\simeq\; \hom(\ounit_{\cC}[\Omega_aA],X),
        \end{equation*}
        and the constant functions and evaluation at the identity loop exhibit $X$ as a retract of $a^*a_*X$, naturally in $X$. Thus, $a_*$ is conservative.
    \end{proof}

    \begin{proposition}
    \label{prop:affineness-and-categorification}
        Let $\cC \in \CAlg(\Pr)$ and $A$ be a truncated anima. Assume that the path animae $\Path_{ab} A$ are $\cC$-affine for all $a,b \in A$. Then $A$ is $\ModC$-affine.
    \end{proposition}
    \begin{proof}
        As $A$ is $\ModC$-ambidextrous~(\cref{exm:ModB-ambi}), by~\cite[Theorem~2.38(3)]{BCSY-Fourier} it is enough to prove that the pushout of the $\cC^A$-algebras $a^* \colon \cC^A \to \cC$ and $b^* \colon \cC^A \to \cC$ is given by
        \begin{equation*}
            \cC \otimes_{\cC^A} \cC \;\simeq\; \cC^{\{a\}\times_A\{b\}} \;\simeq\; \cC^{\Path_{ab} A} \qin \CAlg(\ModC).
        \end{equation*}
        As all the fibers $\Path_{ab} A$ of $a \colon \pt \to A$ are $\cC$-affine,~\cite[Proposition~2.18(3)]{BCSY-Fourier} shows that $a$ is $\cC$-affine.
        Therefore $\cC \simeq \Mod_{a_* \ounit_{\cC}} (\cC^A)$, and the map $a^* \colon \cC^A \to \cC \simeq \Mod_{a_* \ounit_{\cC}}(\cC^A)$ is identified with base change. Therefore
        \begin{equation*}
            \cC \otimes_{\cC^A} \cC \;\simeq\; \Mod_{b^* a_* \ounit_{\cC}} (\cC) \;\simeq\; \Mod_{\ounit_{\cC}^{\Path_{ab} A}}(\cC) \;\simeq\; \cC^{\Path_{ab} A},
        \end{equation*}
        where the second isomorphism is by Beck--Chevalley and the last by affineness.
    \end{proof}

\subsection{Galois descent}
\label{subsec:Galois-descent}
    Fix $\cC \in \CAlg(\Pr)$. 
    \begin{definition}
        Let $R \in \CAlg(\cC^A) \simeq \CAlg(\cC)^A$. Define the local system $\Mod_R(\cC) \in \CAlgC^A$ as the composition
        \begin{equation*}
            A \;\xtoo{R}\; \CAlg(\cC) \;\xtoo{\Mod_{(-)}}\; \CAlgC.
        \end{equation*}
    \end{definition}
    \begin{remark}
        This local system satisfies $a^* \Mod_R(\cC) \simeq \Mod_{R_a}(\cC)$ for any $a \in A$, and the induced symmetric monoidal functor
        \begin{equation*}
            \cC \;\simeq\; a^* A^* \cC \;\too\; a^* \Mod_R(\cC) \;\simeq\; \Mod_{R_a}(\cC)
        \end{equation*}
        is the base-change functor given by tensoring with $R_a$.
    \end{remark}
 
    When $A = \BG$, this construction gives a natural $G$-action on $\Mod_R(\cC)$. Galois descent is then the claim that $\Mod_R(\cC)^{hG} \simeq \cC$.
    \begin{remark}
        The map $\cC \to \Mod_R(\cC)^{hG}$ is the mate of the unit map of $\Mod_R(\cC) \in \CAlgC^{\BG}$. Thus, Galois descent is the first Galois condition~\labelcref{cond:Galois-1} for $\Mod_R(\cC)$.
    \end{remark}
    We fix a \CO-dualizable anima $A \in \An$ and $R \in \Gal_A(\cC)$. We write $F_R \colon \cC^A \to \cC \in \CAlgC$ for the corresponding $\cC$-linear symmetric monoidal functor.
    
    We first reconstruct the local system $\Mod_R(\cC)$ from the categorical interpretation of the Galois extension $R$.

    The following is a pointwise categorical interpretation of Galois descent:
    \begin{proposition}
    \label{prop:categorical-descent-pointwise}
        For every $a \in A$, the following square is a pushout in $\CAlgC$:
        \begin{equation*}
            \begin{tikzcd}
                {\cC^A} & {\cC} \\
                {\cC} & {\Mod_{R_a}(\cC).}
                \arrow["{F_R}", from=1-1, to=1-2]
                \arrow["{a^*}"', from=1-1, to=2-1]
                \arrow[from=1-2, to=2-2]
                \arrow[from=2-1, to=2-2]
                \arrow["\lrcorner"{anchor=center, pos=0.125, rotate=180}, draw=none, from=2-2, to=1-1]
            \end{tikzcd}
        \end{equation*}
    \end{proposition}
    \begin{proof}
        By~\cref{lem:a*-is-affine}, $a^*$ is $\cC$-affine. Hence $\cC \simeq \Mod_{a_* \ounit_{\cC}}(\cC^A)$, under which $a^*$ identifies with the base-change functor
        \begin{equation*}
            a_* \ounit_{\cC} \otimes - \;\colon\; \cC^A \;\too\; \Mod_{a_* \ounit_{\cC}}(\cC^A).
        \end{equation*}
        Its pushout along $F_R$ is
        \begin{equation*}
            \cC \otimes_{\cC^A} \Mod_{a_* \ounit_{\cC}}(\cC^A) \;\simeq\; \Mod_{F_R(a_* \ounit_{\cC})}(\cC).
        \end{equation*}
        Consider the pullback square
        \begin{equation*}
            \begin{tikzcd}
                \pt & A \\
                A & {A\times A.}
                \arrow["a", from=1-1, to=1-2]
                \arrow["a"', from=1-1, to=2-1]
                \arrow["\lrcorner"{anchor=center, pos=0.125}, draw=none, from=1-1, to=2-2]
                \arrow["{(\id, a)}", from=1-2, to=2-2]
                \arrow["\Delta", from=2-1, to=2-2]
            \end{tikzcd}
        \end{equation*}
        The Beck--Chevalley condition gives
        \begin{equation*}
            (\id, a)^* \ounit_{\cC}^{\Path A} \;=\; (\id, a)^* \Delta_* A^* \ounit_{\cC} \;\simeq\; a_* a^* A^* \ounit_{\cC} \;\simeq\; a_* \ounit_{\cC}.
        \end{equation*}
        Using~\cref{lem:lower-shriek-functoriality},
        \begin{equation*}
            F_R(a_* \ounit_{\cC}) \;\simeq\; F_R((\id,a)^* \ounit_{\cC}^{\Path A}) \;\simeq\; a^* F_R^A(\ounit_{\cC}^{\Path A}) \simeq a^* R \;=\; R_a.
        \end{equation*}
    \end{proof}


    \begin{corollary}
    \label{cor:categorical-descent}
        The following is a pushout square in $\CAlgC^A$:
        \begin{equation*}
            \begin{tikzcd}
                {A^* \cC^A} & {A^* \cC} \\
                {A^* \cC} & {\Mod_R(\cC).}
                \arrow["{A^* F_R}", from=1-1, to=1-2]
                \arrow["{\ev_A}"', from=1-1, to=2-1]
                \arrow[from=1-2, to=2-2]
                \arrow[from=2-1, to=2-2]
                \arrow["\lrcorner"{anchor=center, pos=0.125, rotate=180}, draw=none, from=2-2, to=1-1]
            \end{tikzcd}
        \end{equation*}
    \end{corollary}
    \begin{proof}
        The functors $a^* \colon (\ModC)^A \to \ModC$ are symmetric monoidal, colimit-preserving and jointly conservative. The claim therefore follows from~\cref{lem:trivA-at-a-point,prop:categorical-descent-pointwise}.
    \end{proof}

    \begin{theorem}[Galois descent]
    \label{thm:Galois-descent}
        There is a natural isomorphism 
        \begin{equation*}
            \cC \;\isotoo\; A_* \Mod_R(\cC).
        \end{equation*}
        If moreover $\Path_{ab} A$ is $\cC$-affine for all $a,b \in A$, then $\Mod_R(\cC) \in \Gal_A(\ModC)$.
    \end{theorem}

    \begin{proof}
        Let $F_R \colon \cC^A \to \cC$ be the symmetric monoidal functor corresponding to $R$. Consider the unital lax symmetric monoidal functor
        \begin{equation*}
            G \;\colon\; (\ModC)^A \;\xtoo{A_{\scriptscriptstyle\bigstar}}\; \ModC[\cC^A] \;\xtoo{-\otimes_{\cC^A}\cC}\; \ModC,
        \end{equation*}
        where the second functor is base change along $F_R$.
        Its extension over $A$ is
        \begin{equation*}
            G^A \;\colon\; (\ModC)^{A\times A} 
            \;\xtoo{(\pr_2)_{\scriptscriptstyle\bigstar}}\; (\ModC[\cC^A])^A
            \;\simeq\; \Mod_{A^*\cC^A}(\Pr^A)
            \;\xtoo{-\otimes_{A^*\cC^A}A^*\cC}\; (\ModC)^A.
        \end{equation*}
        Evaluating $(\pr_2)_{\scriptscriptstyle\bigstar}$ at $\cC^{\Path A}$ gives
        \begin{equation*}
            (\pr_2)_{\scriptscriptstyle\bigstar}(\cC^{\Path A}) \;\simeq\; (\pr_2)_*\Delta_*A^*\cC \;\simeq\; A^*\cC.
        \end{equation*}
        Hence, by~\cref{cor:categorical-descent},
        \begin{equation*}
            G^A(\cC^{\Path A}) \;\simeq\; A^*\cC\otimes_{A^*\cC^A}A^*\cC \;\simeq\; \Mod_R(\cC).
        \end{equation*}

        On the other hand, $\cC[\Path A]$ is the canonical coGalois extension in $(\ModC)^A$. Thus, the unitality of $G$ and~\cref{lem:lower-shriek-functoriality} give
        \begin{equation*}
            \cC \;\simeq\; G(A^*\cC) \;\simeq\; G((\pr_1)_!\cC[\Path A]) \;\simeq\; A_!G^A(\cC[\Path A]).
        \end{equation*}
        
        There is a natural equivalence
        \begin{equation*}
            \cC[\Path A] \;\simeq\; \cC^{\Path A}
        \end{equation*}
        of underlying objects in $(\ModC)^{A\times A}$. Thus, the preceding calculation gives
        \begin{equation*}
            A_!\Mod_R(\cC) \;\simeq\; \cC.
        \end{equation*}
        The component at $a\in A$ of this colimit cocone is obtained from $a_!\colon\cC\to\cC^A$ by base change along $F_R$. 
        Since
        \begin{equation*}
            F_R(a_!\ounit_{\cC}) \;\simeq\; F_R((a_*\ounit_{\cC})\dual) \;\simeq\; R_a\dual,
        \end{equation*}
        this component is $R_a\dual\otimes_{R_a}-\colon\Mod_{R_a}(\cC)\to\cC$. Its right adjoint is the free-module functor $R_a\otimes-$. Passing to right adjoints turns the colimit cocone into a limit cone, giving the mate of the unit map $A^*\cC\to\Mod_R(\cC)$. Therefore, as colimits in $\Pr$ are computed as limits of the corresponding diagrams of right adjoints,
        \begin{equation*}
            A_* \Mod_R(\cC) \;\simeq\; A_! \Mod_R(\cC) \;\simeq\; \cC.
        \end{equation*}
        
        Assume now that $\Path_{ab} A$ is $\cC$-affine for every $a,b\in A$. By~\cref{lem:affiness-of-ModR}, it is also $\Mod_{R_x}(\cC)$-affine for all $x \in A$. Hence
        \begin{equation*}
            \Mod_R(\cC)\boxtimes\Mod_R(\cC) \;\simeq\; \Mod_{R\boxtimes R}(\cC) \;\simeq\; \Mod_{R^{\Path A}}(\cC)
            \;\simeq\; \Mod_R(\cC)^{\Path A}.
        \end{equation*}
        This is the second Galois condition, while the first one is the descent isomorphism proved above. 
        It remains to prove faithfulness.
        Using the self-duality of $\Mod_R(\cC)$,
        \begin{equation*}
            A_!(\Mod_R(\cC)\dual) \;\simeq\; A_!G^A(\cC[\Path A]) \;\simeq\; \cC.
        \end{equation*}
        Thus $\Mod_R(\cC)$ and $A_!(\Mod_R(\cC)\dual)$ are dualizable. By~\cref{lem:Y-and-f!Y-dual-projection}, for every $\cD \in \ModC$ the projection map is an isomorphism, and hence
        \begin{equation*}
            A_*(A^*\cD \otimes \Mod_R(\cC)) \;\simeq\; \cD \otimes A_*\Mod_R(\cC) \;\simeq\; \cD.
        \end{equation*}
        Therefore 
        \begin{equation*}
            A_* (A^*(-) \otimes \Mod_R(\cC)) \;\simeq\; \id_{\ModC},
        \end{equation*}
        so the base change functor
        \begin{equation*}
            A^*(-) \otimes \Mod_R(\cC) \;\colon\; \ModC \;\too\; (\ModC)^A
        \end{equation*}
        is conservative. Hence $\Mod_R(\cC)$ is faithful.
    \end{proof}

    \begin{remark}
        If $A$ is $\ModC$-affine, then $\Path_{ab} A$ is $\cC$-affine for every $a,b \in A$. In this case, one can see directly that $\Mod_R(\cC)$ is a faithful $A$-Galois extension, since it corresponds to the composition of symmetric monoidal functors
        \begin{equation*}
            G \;\colon\; (\ModC)^A \;\xisotoo{A_{\scriptscriptstyle\bigstar}}\; \ModC[\cC^A] \;\xtoo{- \otimes_{\cC^A} \cC}\; \ModC,
        \end{equation*}
        where the last base change is along $F_R \colon \cC^A \to \cC$.
    \end{remark}

    \begin{lemma}
        The functor obtained from~\cref{thm:Galois-descent},
        \begin{equation*}
            \cC^A \;\simeq\; A_* A^* \cC \;\too\; A_* \Mod_R(\cC) \;\isogetss\; \cC
        \end{equation*}
        is isomorphic to $F_R$.
    \end{lemma}
    \begin{proof}
        The inverse of the descent isomorphism is its right adjoint, and is therefore given by forgetting the module structure and applying $A_*$. Hence the functor in the statement sends $Y\in\cC^A$ to $A_*(R\otimes Y)$, which is $F_R(Y)$ by~\cref{thm:Galois-corepresented}.
    \end{proof}

    Using our categorical interpretation of descent, we can also understand maps of commutative algebras from Galois extensions. These are similar in nature to maps out of Thom algebras (\cite{Antolin-Barthel-2019-Thom}).
    \begin{proposition}
        Let $R \in \Gal_A(\cC)$, $a \in A$ and $S \in \CAlg(\cC)$. Then a map $R_a \to S$ in $\CAlg(\cC)$ is equivalent to an isomorphism $S \otimes R \simeq \triv_a$ in $\Gal_A(\Mod_S(\cC))$. More formally, there is an isomorphism of animae
        \begin{equation*}
            \Map_{\CAlg(\cC)}(R_a, S) \;\simeq\; \{S \otimes R\} \times_{\Gal_A(\Mod_S(\cC))} \{\triv_a\}.
        \end{equation*}
    \end{proposition}

    \begin{proof}
        By the $\Mod_{(-)} \dashv \ounit_{(-)}$ adjunction,
        \begin{equation*}
            \Map_{\CAlg(\cC)}(R_a, S) \;\simeq\; \MapCo(\Mod_{R_a}(\cC), \Mod_S(\cC)).
        \end{equation*}
        Using the pushout formula of \cref{prop:categorical-descent-pointwise}, this is equivalent to the pullback
        \begin{equation*}
            \pt \times_{\MapCo(\cC^A, \Mod_S(\cC))} \pt.
        \end{equation*}
        Under the corepresentability isomorphism $\MapCo(\cC^A,\Mod_S(\cC)) \simeq \Gal_A(\Mod_S(\cC))$, the two points are $S \otimes R$ and $\triv_a$, respectively.
    \end{proof}

\subsection{Relative Galois theory}
\label{subsec:relative-Galois}
    In this subsection, we extend two classical results. 
    Their analogues for faithful Galois extensions of ring spectra indexed by connected 1-finite animae were proved by Rognes~\cite[Theorem~7.2.3]{Rognes-2008-Galois}.
    
    Let $F$ be a field and $E / F$ a finite Galois extension with Galois group $G$. Then if $G \onto H$ is a surjection with kernel $N$, the fixed field $E^{hN} / F$ is a $\BH$-Galois extension. 
    Equivalently, for the map $f \colon \BG \to \BH$, the pushforward $f_* E$ is a $\BH$-Galois extension.
    In the language of rings, this claim remains true even when $G \to H$ is not a surjection, with $H$ finite. In this case, let $N$ be the kernel of the map, so $G/N \subseteq H$, and let $X \coloneq H / (G/N)$. Then $f_* E \simeq \prod_X(E^{hN})$ is a (non-connected) $\BH$-Galois extension.

    Second, if $K \subseteq G$ is a subgroup, then $E / E^{hK}$ is a $\BK$-Galois extension.

    The extension of the first claim is straightforward:
    \begin{proposition}
    \label{thm:relative-quotient}
        Let $f \colon A \to B$ be a map of \CO-dualizable animae.
        Let $R\in \Gal_A(\cC)$. Then $f_*R \in \Gal_B(\cC)$.
    \end{proposition}
    \begin{proof}
        Let $F_R \colon \cC^A \to \cC$ be the $\cC$-linear symmetric monoidal functor corresponding to $R$ and consider the composition of $\cC$-linear symmetric monoidal functors
        \begin{equation*}
            G \;\colon\; \cC^B \;\xtoo{f^*}\; \cC^A \;\xtoo{F_R}\; \cC.
        \end{equation*}
        As $B$ is \CO-dualizable, this represents a $B$-Galois extension given by
        \begin{equation*}
            G^B(\ounit_{\cC}^{\Path B}) \;\simeq\; F_R^B((f \times \id_B)^* (\Delta_B)_* B^* \ounit_{\cC}).
        \end{equation*}
        Consider the pullback square of animae
        \begin{equation*}
            \begin{tikzcd}
                A && B \\
                \\
                {A \times B} && {B \times B.}
                \arrow["f", from=1-1, to=1-3]
                \arrow["{(\id_A,f)}", from=1-1, to=3-1]
                \arrow["\lrcorner"{anchor=center, pos=0.125}, draw=none, from=1-1, to=3-3]
                \arrow["{\Delta_B}", from=1-3, to=3-3]
                \arrow["{f \times \id_B}", from=3-1, to=3-3]
            \end{tikzcd}
        \end{equation*}
        By the Beck--Chevalley condition,
        \begin{equation*}
            (f \times \id_B)^* (\Delta_B)_* B^* \ounit_{\cC} \;\simeq\; (\id_A, f)_* f^* B^* \ounit_{\cC} \;\simeq\; (\id_A \times f)_* (\Delta_A)_* A^* \ounit_{\cC} \;\simeq\; (\id_A \times f)_* \ounit_{\cC}^{\Path A}.
        \end{equation*}
        Similarly, Beck--Chevalley also gives
        \begin{equation*}
            (f \times \id_B)^*\ounit_{\cC}[\Path B] \;\simeq\; (\id_A \times f)_!\ounit_{\cC}[\Path A],
        \end{equation*}
        so $(\id_A \times f)_! \ounit_{\cC}[\Path A]$ is dualizable, with dual
        \begin{equation*}
            ((\id_A \times f)_! \ounit_{\cC}[\Path A])\dual \;\simeq\; (\id_A \times f)_* \ounit_{\cC}^{\Path A}.
        \end{equation*}
        By~\cref{lem:lower-shriek-functoriality},
        \begin{equation*}
            G^B(\ounit_{\cC}^{\Path B}) 
            \;\simeq\; F_R^B((\id_A \times f)_! \ounit_{\cC}[\Path A])\dual
            \;\simeq\; (f_! F_R^A(\ounit_{\cC}[\Path A]))\dual
            \;\simeq\; f_* R.
        \end{equation*}
        Here the last isomorphism uses~\cref{cor:internal-adjoints-maps-of-animae} and $F_R^A(\ounit_{\cC}[\Path A])\dual \simeq F_R^A(\ounit_{\cC}^{\Path A}) \simeq R$. These isomorphisms respect the commutative algebra structures, and therefore identify the represented extension with $f_*R$.
    \end{proof}

    \begin{remark}
    \label{rmrk:functoriality-of-galois}
        In particular, the restriction of the functor
        \begin{equation*}
            \MapCo(\cC^{(-)}, \cC) \;\colon\; \An \;\too\; \An
        \end{equation*}
        with its $(-)^*$-functoriality to \CO-dualizable animae gives $\Gal_{(-)}(\cC)$ the structure of a functor, where the functoriality is given by $(-)_*$.
    \end{remark}

    We now move to the second goal of this subsection, extending the fact that if $E / F$ is $\BG$-Galois and $K \subseteq G$, then $E / E^{hK}$ is $\BK$-Galois. Here the role of a subgroup is replaced by an affine map. This is a relative version of the Galois consequence of affineness in~\cite[Corollary~2.28]{BCSY-Fourier}.
    \begin{proposition}
    \label{thm:relative-subgroup}
        Let $g \colon B \to A$ be a $\cC$-affine map of \CO-dualizable animae.
        Let $R\in\Gal_A(\cC)$. Then $g^* R$ extends naturally to a $B$-Galois extension $g^* R \in \Gal_B(\Mod_{B_* g^* R}(\cC))$.
    \end{proposition}
    \begin{proof}
        As $g$ is $\cC$-affine, we have an isomorphism
        \begin{equation*}
            g_{\scriptscriptstyle\bigstar} \;\colon\; \cC^B \;\isotoo\; \Mod_{g_* B^* \ounit_{\cC}}(\cC^A).
        \end{equation*}
        Composing with the functor on modules induced by $F_R \colon \cC^A \to \cC$, we get the $\cC$-linear symmetric monoidal functor
        \begin{equation*}
            G \;\colon\; \cC^B \;\xisotoo{g_{\scriptscriptstyle\bigstar}}\; \Mod_{g_* B^* \ounit_{\cC}}(\cC^A) \;\xtoo{F_R}\; \Mod_{F_R(g_* B^* \ounit_{\cC})}(\cC).
        \end{equation*}
        
        Using the formula for $F_R$ from~\cref{thm:Galois-corepresented} and the right projection formula (\cref{lem:projection-for-dbl,prop:Gal-dualizable}),
        \begin{equation*}
            F_R(g_* B^* \ounit_{\cC}) \;\simeq\; A_*(R \otimes g_* B^* \ounit_{\cC}) \;\simeq\; A_*g_*g^*R \;\simeq\; B_*g^*R \qin \CAlg(\cC).
        \end{equation*}
        The functor $G$ is then a $\cC$-linear symmetric monoidal functor
        \begin{equation*}
            G \;\colon\; \cC^B \;\xisotoo{g_{\scriptscriptstyle\bigstar}}\; \Mod_{g_* B^* \ounit_{\cC}}(\cC^A) \;\xtoo{F_R}\; \Mod_{B_* g^* R}(\cC).
        \end{equation*}
        As $B$ is \CO-dualizable, this represents a $B$-Galois extension in $\Mod_{B_* g^* R}(\cC)$ given by $G^B(\ounit_{\cC}^{\Path B})$.
        $G^B$ is described as the composition
        \begin{equation*}
            G^B \;\colon\; \cC^{B \times B} \;\xtoo{(g \times \id_B)_{\scriptscriptstyle\bigstar}}\; \Mod_{(g \times \id_B)_* (B \times B)^* \ounit_{\cC}}(\cC^{A \times B}) \;\simeq\; \Mod_{g_* B^* \ounit_{\cC}}(\cC^A)^B \;\xtoo{F_R^B}\; \Mod_{B_* g^* R}(\cC)^B. 
        \end{equation*}
        Recall that $(g \times \id_B)_{\scriptscriptstyle\bigstar}$ is computed as $(g \times \id_B)_*$ remembering the $(g \times \id_B)_* (B \times B)^* \ounit_{\cC}$-module structure. Therefore
        \begin{equation*}
            (g \times \id_B)_{\scriptscriptstyle\bigstar}(\ounit_{\cC}^{\Path B}) \;\simeq\; (g \times \id_B)_* (\Delta_B)_* B^* \ounit_{\cC} \;\simeq\; (g, \id_B)_* B^* \ounit_{\cC}.
        \end{equation*}
        Consider the pullback square
        \begin{equation*}
            \begin{tikzcd}
                B && {A \times B} \\
                \\
                A && {A\times A.}
                \arrow["{(g,\id_B)}", from=1-1, to=1-3]
                \arrow["g", from=1-1, to=3-1]
                \arrow["\lrcorner"{anchor=center, pos=0.125}, draw=none, from=1-1, to=3-3]
                \arrow["\id_A \times g", from=1-3, to=3-3]
                \arrow["{\Delta_A}", from=3-1, to=3-3]
            \end{tikzcd}
        \end{equation*}
        By Beck--Chevalley and~\cref{lem:lower-shriek-functoriality},
        \begin{equation*}
            \begin{split}
                G^B(\ounit_{\cC}^{\Path B})
                & \;\simeq\; F_R^B((g, \id_B)_* B^* \ounit_{\cC}) \\
                & \;\simeq\; F_R^B((g, \id_B)_* g^* A^* \ounit_{\cC}) \\
                & \;\simeq\; F_R^B((\id_A \times g)^* (\Delta_A)_* A^* \ounit_{\cC}) \\
                & \;\simeq\; g^* F_R^A(\ounit_{\cC}^{\Path A}) \simeq g^* R.
            \end{split}
        \end{equation*}
        Thus $g^* R \in \Gal_B(\Mod_{B_* g^* R}(\cC))$.
    \end{proof}

\section{Galois prorepresentability}
\label{sec:Galois-anima}
Throughout this section, we fix a small full subcategory $\cT\subseteq \An$ which is closed under finite limits and consists of truncated, \COB-dualizable and $\ModB$-affine animae.
We show in~\cref{subsec:prorepresentability} that under these assumptions there exists a pro-$\cT$ object $\piT(\cC)$ representing $\cT$-Galois extensions. 
In~\cref{subsec:categorification} we study its behavior under categorification. 
In~\cref{subsec:pro-Galois} we define pro-Galois extensions and extend the representability result to this context.

\subsection{The Galois pro-anima}
\label{subsec:prorepresentability}
 
    \begin{theorem}
    \label{thm:pro-affine-anima}
        There exists a functor
        \begin{equation*}
            \piT \;\colon\; \CAlgB \;\too\; \Pro(\cT)\op
        \end{equation*}
        prorepresenting Galois extensions, i.e.\ for any $\cC \in \CAlgB$ and $A \in \cT$, there is a natural isomorphism of animae
        \begin{equation*}
            \Map(\piT(\cC), A) \;\simeq\; \Gal_A(\cC).
        \end{equation*}
    \end{theorem}
    
    \begin{proof}
        By~\cref{thm:Galois-corepresented,rmrk:functoriality-of-galois}, there exists a functor
        \begin{equation*}
            \begin{split}
                \Gal_{(-)}(-) \;\colon\; \CAlgB & \;\too\; \Fun(\cT, \An) \\
                \cC & \;\mapsto\; (A \mapsto \Gal_A(\cC)),
            \end{split}
        \end{equation*}
        which is identified with the functor
        \begin{equation*}
            \begin{split}
                \MapBo(\basecat^{(-)}, -) \;\colon\; \CAlgB & \;\too\; \Fun(\cT, \An) \\
                \cC & \;\mapsto\; (A \mapsto \MapBo(\basecat^A, \cC)).
            \end{split}
        \end{equation*}
        To show that it lands in $\Pro(\cT)\op$, we must show that for any $\cC \in \CAlgB$ the corresponding functor preserves finite limits (\cite[Corollary~5.3.5.4]{Lurie-HTT}), or equivalently, sends the terminal object $\pt \in \cT$ to the point and preserves pullbacks. The first is clear. Let $A \to B \leftarrow C$ be a pullback diagram in $\cT$. Since $B$ and the fibers of $A \to B$ are truncated, $B$ and the map are $\ModB$-ambidextrous (\cref{exm:ModB-ambi}), while $B$ is $\ModB$-affine. Thus,~\cite[Theorem~2.38]{BCSY-Fourier} shows that $A \to B$ is Eilenberg--Moore. That is
        \begin{equation*}
            \basecat^A \otimes_{\basecat^B} \basecat^C \;\simeq\; \basecat^{A \times_B C}.
        \end{equation*}
        Applying $\MapBo(-,\cC)$ gives the required pullback.
        Therefore it lands in prorepresentable objects.
    \end{proof}
    
    \begin{notation}
        We call $\piT(\cC)$ the $\cT$-Galois anima of $\cC$.
        For $R \in \Gal_A(\cC)$, we denote the corresponding continuous map by $\chi_R \colon \piT(\cC) \to A$.
    \end{notation}
        
    \begin{example}
        Let $A \in \cT$ and $a \in A$. The character of the trivial Galois extension $\triv_a$ is the trivial character which we also denote
        \begin{equation*}
            \triv_a \;\colon\; \piT(\cC) \;\too\; \pt \;\xtoo{a}\; A.
        \end{equation*}
    \end{example}
    
    \begin{example}
    \label{exm:galois-anima-of-local-systems}
        For every $A \in \cT$, there is an isomorphism $\piT(\cC^A) \simeq \piT(\cC)[A] = \colim_A \piT(\cC)$.
    \end{example}
    \begin{proof}
        Let $B \in \cT$. Then
        \begin{equation*}
            \MapBo(\basecat^B, \cC^A) \;\simeq\; \MapBo(\basecat^B, \cC)^A \;\simeq\; \Map(\piT(\cC), B)^A \;\simeq\; \Map(\piT(\cC)[A], B).
        \end{equation*}
    \end{proof}
    
    \begin{example}[Pointed categories]
        Let $\cC \in \CAlg(\Pr)$ be pointed. Then all finite sets are \CO-dualizable. All finite sets are also $\ModC$-affine by~\cite[Example~2.35]{BCSY-Fourier}.
        Therefore, there exists $\pim[0](\cC) \in \Pro(\Set^{\mrm{fin}})$ representing Galois extensions along finite sets.
    \end{example}

    \begin{example}
    \label{exm:classical-case}
        Let $\cC$ be a semiadditive category. 
        Then all finite sets are dualizable in $\cC$, hence all $1$-finite animae are \CO-dualizable. 
        Finite sets are $\cC$-affine by~\cite[Example~2.35]{BCSY-Fourier}. Since the path animae of a $1$-finite anima are finite sets, \cref{prop:affineness-and-categorification} shows that all $1$-finite animae are also $\ModC$-affine.
        Therefore, there exists a pro-$1$-finite Galois anima $\pim[1] \cC$ classifying Galois extensions.        
        When $\cC = \Mod_R(\Ab)$, it agrees with the \'etale fundamental groupoid $\pi^{\et}_{\le 1} (\Spec R)$. 
        When $\cC$ is stable, it agrees with Mathew's profinite groupoid $\pi_{\le 1}^{\mrm{weak}} \cC$~(\cite[Definition~6.8]{Mathew-2016-Galois}), since by~\cite[Proposition~6.13]{Mathew-2016-Galois} they both satisfy the same universal property.
    \end{example}
    
    Let $\cT' \subseteq \cT$ be a subcategory closed under finite limits. Then the inclusion $\Pro(\cT') \to \Pro(\cT)$ preserves all limits. By the adjoint functor theorem (for the presentable categories $\Pro(\cT)\op = \Ind(\cT\op)$ and $\Pro(\cT')\op = \Ind(\cT'{}\op)$) it admits a left adjoint.
    \begin{definition}
    \label{def:protruncation-of-Ts}
        Let $\cT' \subseteq \cT$ be a subcategory closed under finite limits. Define the truncation functor
        \begin{equation*}
            \tau_{\scriptscriptstyle\cT'} \;\colon\; \Pro(\cT) \to \Pro(\cT')
        \end{equation*}
         as the left adjoint of the inclusion $\Pro(\cT') \to \Pro(\cT)$.
    \end{definition}
    
    \begin{lemma}
    \label{lem:T-truncation-of-Galois-anima}
        Let $\cT' \subseteq \cT$ be a subcategory closed under finite limits. Then for any $\cC \in \CAlgB$
        \begin{equation*}
            \tau_{\scriptscriptstyle\cT'} \piT(\cC) \;\simeq\; \piT[\cT'](\cC).
        \end{equation*}
    \end{lemma}
    \begin{proof}
        Let $A \in \cT'$. Then
        \begin{equation*}
            \Map_{\Pro(\cT')}(\piT[\cT'](\cC), A) \;\simeq\; \Gal_A(\cC) \;\simeq\; \Map_{\Pro(\cT)}(\piT(\cC), A) \;\simeq\; \Map_{\Pro(\cT')}(\tau_{\scriptscriptstyle\cT'} \piT(\cC), A).
        \end{equation*}
    \end{proof}

\subsection{Categorification}
\label{subsec:categorification}
    We show that categorifications of $\cC \in \CAlg(\ModB)$ admit a pro-Galois anima extending that of $\cC$, which classifies a larger class of Galois extensions.

    We use the notion of presentable $2$-categories as developed in~\cite{stefanich-2020-presentable,Aoki-2025-higher-presentable,Scholze-2026-Gestalten}.
    We assume now $\Pr = \PrL_{\kappa}$ for some uncountable regular cardinal $\kappa$, so that $\Pr \in \CAlg(\Pr)$~(\cite{Campion-2022-presentable},~\cite[Proposition~2.3]{Aoki-2025-higher-presentable}, see also~\cite[Theorem~3.11]{Scholze-2026-Gestalten}). We then define $2\Pr \coloneq \ModC[\Pr] = \Mod_{\Pr}(\Pr)$.
    There is an adjunction
    \begin{equation*}
        \ModC[(-)] \;\colon\; \CAlg(\Pr) \;\longadj\; \CAlg(2\Pr) \;\cocolon\; \Omega,
    \end{equation*}
    where for a presentable $2$-category $\cD$, its loops $\Omega \cD \coloneq \End(\ounit_{\cD})$ is the $1$-category of endomorphisms of the unit, considered as an object of $\CAlg(\Pr)$.
    
    \begin{definition}
        Let $\cC \in \CAlg(\Pr)$. A categorification of $\cC$ is a presentably symmetric monoidal $2$-category $\cD \in \CAlg(2\Pr)$ with an isomorphism $\Omega \cD \simeq \cC$.
    \end{definition}
    
    The universal categorification is $\cD = \ModC$. The isomorphism $\cC \isoto \Omega \cD$ admits a mate $u \colon \ModC \to \cD$, showing that $\cD \in 2\ModC \coloneq \Mod_{\ModC}(\Pr)$. In particular
    \begin{corollary}
        If $\cC \in \CAlgB$ and $\cD$ is a categorification of $\cC$, then $\cD \in \CAlg(2\ModB)$.
    \end{corollary}

    \begin{definition}
        Choose an infinite regular cardinal $\lambda$ such that $|\pi_0 A| < \lambda$ for every $A\in\cT$.
        Let $\cT^+_{\lambda} \subseteq \An$ be the full subcategory of truncated animae $A$ with $|\pi_0 A|<\lambda$ such that the path animae $\Path_{ab} A \in \cT$ for all $a,b \in A$.
    \end{definition}
    In particular, as $\cT$ is closed under finite limits, $\cT \subseteq \cT^+_{\lambda}$, and $\cT^+_{\lambda}$ is closed under finite limits.
    The bound on the number of components keeps $\cT^+_{\lambda}$ a small category.
    
    \begin{lemma}
    \label{lem:T-is-also-for-ModC}
        Every anima in $\cT^+_{\lambda}$ is \CO[\ModB]-dualizable and $2\ModB = \ModC[\ModB]$-affine.
    \end{lemma}
    \begin{proof}
        By~\cref{exm:ModB-ambi}, all truncated animae are $\ModB$-ambidextrous. In particular all animae in $\cT^+_{\lambda}$ are \CO[\ModB]-dualizable.
        Moreover, for every $A \in \cT^+_{\lambda}$ and $a,b \in A$, the path anima $\Path_{ab} A$ lies in $\cT$ and is therefore $\ModB$-affine. Thus, \cref{prop:affineness-and-categorification} shows that $A$ is $2\ModB$-affine.
    \end{proof}
    
    \begin{corollary}
        Let $\cD$ be a categorification of $\cC$.
        There exists a pro-anima $\piT[\cT^+_{\lambda}](\cD) \in \Pro(\cT^+_{\lambda})$ classifying $\cT^+_{\lambda}$-Galois extensions in $\cD$.
    \end{corollary}
    \begin{proof}
        This follows from~\cref{thm:pro-affine-anima} and~\cref{lem:T-is-also-for-ModC}, using that $\cD \in \CAlg(2\ModB)$.
    \end{proof}

    \begin{remark}
        Define $n\Pr^{\oplus}$ inductively as $\Mod_{(n-1)\Pr^{\oplus}}(\Pr)$, where $0\Pr^{\oplus} = \CMon$. Thus $1\Pr^{\oplus}$ is the category of semiadditive presentable categories.
        Then, choosing iteratively $\lambda = \omega_0$, any category $\cC \in \CAlg(n\Pr^{\oplus})$ admits an $n$-truncated pro-$\pi$-finite anima $\pim[n]\cC$ classifying $n$-truncated $\pi$-finite Galois extensions in $\cC$.
    \end{remark}
    
    \begin{proposition}
    \label{prop:piT-and-categorification}
        Let $\cD$ be a categorification of $\cC$.
        There is an isomorphism $\tau_{\scriptscriptstyle\cT}{\piT[\cT^+_{\lambda}]}(\cD) \simeq \piT(\cC)$. In particular, for any $A \in \cT$ there is an isomorphism
        \begin{equation*}
            \begin{split}
                \Gal_A(\cC) & \;\isotoo\; \Gal_A(\cD) \\
                R & \;\mapsto\; u(\Mod_R(\cC)),
            \end{split}
        \end{equation*}
        where $u \colon \ModC \to \cD$ is the natural map.
    \end{proposition}
    \begin{proof}
        By~\cref{lem:T-truncation-of-Galois-anima} it is equivalent to showing that $\piT(\cD) \simeq \piT(\cC)$.
        Consider the functor $\Gal_{(-)}(\cD) \colon \cT \to \An$ given by
        \begin{equation*}
            \Gal_{(-)}(\cD) = \MapCo[\ModC]((\ModC)^{(-)}, \cD) \;\colon\; \cT \;\too\; \An.
        \end{equation*}
        All objects in $\cT$ are $\ModB$-affine and hence, by~\cref{lem:affiness-of-ModR}, they are also $\Mod_{\cC}(\ModB) = \ModC$-affine. Therefore there is a natural isomorphism
        \begin{equation*}
            (\ModC)^{(-)} \;\isoToo\; \ModC[\cC^{(-)}] \qin \Fun(\cT\op, \CAlgC[\ModC]).
        \end{equation*}
        Using the symmetric monoidal adjunction $\ModC[(-)] \dashv \Omega$, we obtain
        \begin{equation*}
            \Gal_{(-)}(\cD) \;=\; \MapCo[\ModC]((\ModC)^{(-)}, \cD) \;\simeq\; \MapCo[\ModC](\ModC[\cC^{(-)}], \cD) \;\simeq\; \MapCo(\cC^{(-)}, \cC) \;\simeq\; \Gal_{(-)}(\cC).
        \end{equation*}
        This in particular implies $\piT(\cD) \simeq \piT(\cC)$.
    
        Moreover, for any $A \in \cT$, the isomorphism $\Gal_A(\cC) \isoto \Gal_A(\cD)$ is the composite
        \begin{equation*}
            \begin{split}
                \Gal_A(\cC) 
                & \;\isotoo\; \MapCo(\cC^A, \cC) \\
                & \;\xisotoo{u \circ \Mod_{(-)}}\; \MapCo[\ModC](\ModC[\cC^A], \cD) \\
                & \;\xisotoo{\circ A_{\scriptscriptstyle\bigstar}}\; \MapCo[\ModC]((\ModC)^A, \cD) \\
                & \;\simeq\; \Gal_A(\cD).
            \end{split}
        \end{equation*}
        The identification in the proof of~\cref{thm:Galois-descent} sends $R \in \Gal_A(\cC)$ to $\Mod_R(\cC) \in \Gal_A(\ModC)$, and then to $u(\Mod_R(\cC)) \in \Gal_A(\cD)$ by composition.
    \end{proof}

\subsection{Pro-Galois extensions}
\label{subsec:pro-Galois}
    We now extend the previous results to pro-objects. We write an object of $\Pro(\cT)$ as a formal cofiltered limit $A = \qlimi A_i$.

    \begin{definition}
    \label{def:continuous-local-systems}
        For $A \in \cT$, we equip $\basecat^A$ with the pointwise symmetric monoidal structure. We then extend the functor 
        \begin{equation*}
            \basecat^{(-)} \;\colon\; \cT\op \;\too\; \CAlgB
        \end{equation*}
        to a filtered-colimit-preserving functor
        \begin{equation*}
            \basecat_{\cts}^{(-)} \;\colon\; \Pro(\cT)\op \;\too\; \CAlgB.
        \end{equation*}
        It sends $A = \qlimi A_i$ to $\basecat_{\cts}^A = \colim_i \basecat^{A_i}$, computed in $\CAlgB$. 
        More generally, for $\cC \in \CAlgB$, we write $\cC_{\cts}^A \coloneq \colim_i \cC^{A_i}$. We refer to objects of $\cC_{\cts}^A$ as \emph{continuous $A$-local systems}.

        One can avoid a choice of a presentation of $A$ by choosing the canonical index category $\cT\op \times_{\Pro(\cT)\op} (\Pro(\cT)_{A/})\op$ (\cite[Corollary~5.3.5.4]{Lurie-HTT}).
    \end{definition}

    The same construction in $\Pr$ defines $\cD_{\cts}^A$ for any presentable category $\cD$ and any pro-anima $A$.


    \begin{remark}
    \label{rmk:filtered-colimit-limit}
        Since the colimit $\cC_{\cts}^A = \colim_i \cC^{A_i}$ is filtered, it is computed in the underlying category $\ModB$ by~\cite[Corollary~3.2.3.2]{Lurie-HA}. As a category, $\cC_{\cts}^A$ is equivalently the limit $\lim_i \cC^{A_i}$ along the right adjoints $(p_{ij})_*$, where $p_{ij} \colon A_i \to A_j$ are the transition maps. Concretely, an object $X \in \cC_{\cts}^A$ is a compatible system $(X_i \in \cC^{A_i})_i$ with isomorphisms $(p_{ij})_* X_i \simeq X_j$ for every transition map $p_{ij} \colon A_i \to A_j$, and all the required coherences.
    \end{remark}

    \begin{definition}
        Let $f \colon A \to B$ be a map of pro-animae. By the functoriality of $\basecat_{\cts}^{(-)}$ there exists a functor
        \begin{equation*}
            f^* \;\colon\; \basecat_{\cts}^B \;\too\; \basecat_{\cts}^A \qin \CAlgB.
        \end{equation*}
        Denote its lax monoidal right adjoint by $f_*$.
    \end{definition}
    \begin{remark}
        After a common cofiltered reindexing, write $f$ as a compatible system of maps $f_i\colon A_i\to B_i$. Denote the projections by $p_i\colon A\to A_i$ and $q_i\colon B\to B_i$. For $X=(X_i)\in\cC_{\cts}^A$, the pushforward is the system
        \begin{equation*}
            f_*X \;\simeq\; \big((f_i)_*X_i\big)_i \qin \cC_{\cts}^B.
        \end{equation*}

        For $Y=(Y_i)\in\cC_{\cts}^B$, the pullback is
        \begin{equation*}
            f^*Y \;\simeq\; \colim_i p_i^*f_i^*Y_i \qin \cC_{\cts}^A.
        \end{equation*}
        Indeed, $Y\simeq\colim_i q_i^*Y_i$ (see~\cref{rmrk:ind-local-system}). Apply $f^*$ and use $f^*q_i^*\simeq p_i^*f_i^*$.
    \end{remark}

    \begin{definition}
        Define the functor
        \begin{equation*}
            \ProGal_{(-)}(\cC) \;\colon\; \Pro(\cT) \;\too\; \An
        \end{equation*} 
        to be the unique cofiltered-limit-preserving functor extending $\Gal_{(-)}(\cC) \colon \cT \to \An$. Thus, $\ProGal_A(\cC) = \lim_i \Gal_{A_i}(\cC)$ for $A = \qlimi A_i$.
    
        An object $R \in \ProGal_A(\cC)$ is called an \emph{$A$-pro-Galois extension}. The inclusion 
        \begin{equation*}
            \ProGal_A(\cC) \;=\; \lim_i\Gal_{A_i}(\cC) \;\subseteq\; \lim_i \CAlg(\cC)^{A_i} \;=\; \CAlg(\cC)^A_{\cts}
        \end{equation*}
        identifies it as a compatible system of $R_i \in \Gal_{A_i}(\cC)$ with isomorphisms $(p_{ij})_* R_i \simeq R_j$ for every transition map $p_{ij} \colon A_i \to A_j$.
    \end{definition}

    The functoriality of $\ProGal_{(-)}(\cC)$ is the right Kan extension of the $(-)_*$-functoriality of $\Gal_{(-)}(\cC)$. It therefore agrees with the $(-)_*$-functoriality of continuous local systems.

    \begin{corollary}
        Let $f \colon A \to B \in \Pro(\cT)$ and $R \in \ProGal_A(\cC)$. Then $f_* R \in \ProGal_B(\cC)$.
    \end{corollary}

    \begin{proposition}
    \label{prop:pro-Galois-classification}
        There are isomorphisms, natural in $A \in \Pro(\cT)$ and $\cC \in \CAlgB$,
        \begin{equation*}
                \ProGal_A(\cC)
                \;\simeq\; \MapBo\big(\basecat_{\cts}^A, \cC\big)
                \;\simeq\; \Map_{\Pro(\cT)}\big(\piT(\cC), A\big).
        \end{equation*}
        In particular, $\piT(\cC)$ represents pro-Galois extensions.
    \end{proposition}
    \begin{proof}
        All three functors preserve cofiltered limits in $A$ and agree on $\cT$ by~\cref{thm:Galois-corepresented,thm:pro-affine-anima}.
    \end{proof}

    For $A \in \cT \subseteq \Pro(\cT)$ one has $\ProGal_A(\cC) \simeq \Gal_A(\cC)$ and $\basecat^A_{\cts} \simeq \basecat^A$, recovering~\cref{thm:Galois-corepresented,thm:pro-affine-anima}.

    \begin{corollary}
    \label{cor:piT-right-adjoint}
        The functor $\piT \colon \CAlgB \to \Pro(\cT)\op$ is right adjoint to $\basecat_{\cts}^{(-)}$.
    \end{corollary}

    \begin{corollary}
       For $\cC, \cD \in \CAlgB$, there is an isomorphism $\piT(\cC \times \cD) \simeq \piT(\cC) \sqcup \piT(\cD)$.
    \end{corollary}

    \begin{remark}
    \label{rmrk:ind-local-system}
        Let $\LocSys(\cC)$ be the cocartesian Grothendieck construction of $\cC^{(-)} \colon \An\op \to \Pr$. Its objects are pairs $(A, X)$ where $A \in \An$ and $X \in \cC^A$. 
        Similarly, let $\LocSys_{\cts}(\cC)$ be the cocartesian Grothendieck construction of $\cC_{\cts}^{(-)} \;\colon (\ProAn)\op \to \Pr$. Its objects are $(A, X)$ where $A \in \ProAn$ and $X \in \cC^A_{\cts}$.

        The inclusion $\An \into \ProAn$ induces a functor $\LocSys(\cC) \into \LocSys_{\cts}(\cC)$. Formally extending it, we get a functor
        \begin{equation*}
            \Ind\LocSys(\cC) \;\too\; \LocSys_{\cts}(\cC)
        \end{equation*}
        sending $\qcolimi(A_i,X_i)$ to
        \begin{equation*}
            \big(A,\colim_i p_i^*X_i\big), \qquad A=\qlimi A_i,
        \end{equation*}
        where $p_i\colon A\to A_i$ are the projections and $p_i^*\colon\cC^{A_i}\to\cC_{\cts}^A$ are the canonical functors. This functor admits a right adjoint
        \begin{equation*}
            \LocSys_{\cts}(\cC) \;\too\; \Ind\LocSys(\cC)
        \end{equation*}
        that for $A = \qlimi A_i$ and $X = (X_i) \in \cC^A_{\cts}$ sends the pair $(A,X)$ to $\qcolimi (A_i, X_i)$.

        For a compatible system $X=(X_i)$, the counit gives an isomorphism
        \begin{equation*}
            \colim_i p_i^*X_i \;\isotoo\; X.
        \end{equation*}
        Therefore $\LocSys_{\cts}(\cC)$ embeds fully faithfully into $\Ind\LocSys(\cC)$.

    \end{remark}

    The limit, or materialization~(\cite[Remark~E.2.1.4]{Lurie-SAG}) functor
    \begin{equation*}
        \lim \;\colon\; \Pro(\cT) \;\too\; \An
    \end{equation*}
    is the unique limit-preserving functor extending the embedding $\cT \to \An$. It sends $A = \qlimi A_i$ to the actual limit $\lim_i A_i$ computed in $\An$.

    \begin{remark}
    \label{rmk:galois-vs-progalois-clean}
        Any $A$-pro-Galois extension has an underlying $\lim A$-local system of algebras.
        Write $A = \qlimi A_i \in \Pro(\cT)$ and $R \in \ProGal_A(\cC)$ with system $(R_i)$.
        Its \emph{underlying local system} is
        \begin{equation*}
            \underline{R} \;\coloneq\; \colim_i p_i^* R_i \;\in\; \CAlg(\cC)^{\lim A},
        \end{equation*}
        where $p_i \colon \lim A \to A_i$ are the projections.
        Indeed, for every transition map $p_{ij} \colon A_i \to A_j$, the counit of $p_{ij}^* \dashv (p_{ij})_*$ gives a map
        \begin{equation*}
            p_j^* R_j \;\simeq\; p_i^*p_{ij}^*(p_{ij})_*R_i \;\too\; p_i^*R_i,
        \end{equation*}
        and these maps form the filtered diagram in the colimit.
        This is the value
        on $R$ of the forgetful functor
        $\ProGal_A(\cC) \to \CAlg(\cC)^{\lim A}$ induced by $\Gal_{(-)}(\cC) \to \CAlg(\cC)^{(-)}$.

        If $A$ is constant, then $A = \lim A$ and $\underline{R} = R$. Otherwise, even when $\lim A \in \cT$, $R$ being $A$-pro-Galois is very different from $\underline{R}$ being $\lim A$-Galois.
    \end{remark}

\bibliographystyle{alpha}
\phantomsection\addcontentsline{toc}{section}{\refname}
\bibliography{references}

\end{document}